\documentclass[leqno]{article}
\usepackage{amssymb,amsmath,amsfonts,amscd,verbatim}

\newtheorem{theorem}{Theorem}[section]
\newtheorem{corollary}[theorem]{Corollary}
\newtheorem{lemma}[theorem]{Lemma}
\newtheorem{proposition}[theorem]{Proposition}
\newtheorem{definition}[theorem]{Definition}
\newtheorem{remark}[theorem]{Remark}

\numberwithin{equation}{section}

\def\sqr#1#2{{\vcenter{\vbox{\hrule height.#2pt
    \hbox{\vrule width.#2pt height#1pt \kern#1pt
    \vrule width.#2pt}
    \hrule height.#2pt}}}}

\def\ga{\gamma}

\def\De{\Delta}
\def\de{\delta}

\def\la{\lambda}
\def\eps{\varepsilon}

\def\om{\omega}
\def\Om{\Omega}

\def\pa{\partial}

\def\bB{\bar{B}}

\def\tc{\tilde c}

\font\bbb=msbm10

\def\R{\hbox{{\bbb R}}}

\def\cL{{\cal L}}

\def\capa{\, {\rm cap}\, }
\def\mes{\,{\rm mes}\, }
\def\Lip{\, {\rm Lip}}

\def\tmu{\tilde{\mu}}

\def\cF{{\mathcal F}}
\def\cG{{\mathcal G}}

\def\ms{\medskip}
\def\bs{\bigskip}

\begin{document}

\title{Discreteness of spectrum and strict positivity
criteria for magnetic Schr\"odinger
operators.}
\author{{\bf Vladimir Kondratiev\thanks{Research partially supported by 
Oberwolfach Forschungsinstitut f\"ur Mathematik and Northeastern University}}
\smallskip\\
Department of Mechanics and Mathematics\\
Moscow State University\\
Vorobievy Gory, Moscow, 119899,
Russia\medskip\\
E-mail: kondrat@vnmok.math.msu.su
\bigskip\\
{\bf Vladimir Maz'ya\thanks{Research partially supported by the Department of
Mathematics and the Robert G. Stone Fund at Northeastern
University}}
\smallskip\\
Department of Mathematics\\
Link\"oping University\\
Sweden\medskip\\
E-mail: vlmaz@mai.liu.se
\bigskip\\
{\bf Mikhail Shubin
\thanks{Research partially supported by NSF grant DMS-0107796} }
\\
Department of Mathematics\\
Northeastern University\\
Boston, MA 02115,
USA
\medskip\\
E-mail: shubin@neu.edu
}

\date{}

\maketitle

\begin{abstract}

We establish necessary and sufficient conditions
for the discreteness of spectrum and strict positivity
of magnetic Schr\"odinger operators with a positive scalar potential.  
They are expressed in terms of Wiener's capacity and 
the local energy of the magnetic field. 
The conditions for the discreteness of spectrum depend, in particular, on  
a functional parameter which is a decreasing function of one variable
whose argument is the normalized local energy of the magnetic field. 
This function  enters  the negligibility condition of sets for the scalar potential.
We give a description for the range of all admissible functions
which is precise in a certain sense. 

In case when there is no magnetic field, our results extend 
the discreteness of spectrum and positivity criteria by A.~Molchanov (1953)
and V.~Maz'ya (1973).
\end{abstract}

\section{Introduction and main results}\label{S:Intro}

The main object of this paper is the magnetic Schr\"odinger operator
in $\R^n$ which has the form

\begin{equation}\label{E:Schroed}
H_{a,V}=\sum_{j=1}^n P_j^2+V,
\end{equation}
where
\begin{equation*}\label{E:Pj}
P_j=\frac{1}{i}\frac{\pa}{\pa x^j}+a_j,
\end{equation*}
and $a_j=a_j(x)$, $V=V(x)$, $x=(x^1,\dots,x^n)\in\R^n$.
We assume that $a_j$ and $V$ are real-valued functions. Denote also
\begin{equation*}\label{E:nabla_a}
\nabla_a u=\nabla u +iau=
\left(\frac{\pa u}{\pa x^1}+ia_1 u,\dots, \frac{\pa u}{\pa x^n}+ia_n u\right).
\end{equation*}

We will assume a priori that $V\in L^1_{loc}(\R^n)$ and $a\in L^2_{loc}(\R^n)$
(which will be a shorthand for saying that $a_j\in L^2_{loc}(\R^n)$ for all $j=1,\dots,n$). 
This allows to define the quadratic form 
\begin{equation}\label{E:haV}
h_{a,V}(u,u)=\int_{\R^n}(|\nabla_a u|^2+V|u|^2)dx
\end{equation}
on functions $u\in C_c^\infty(\R^n)$. A stronger local requirement on $a$ will be imposed
for the discreteness of spectrum results. (For example, it will be sufficient to
require that $a\in L^\infty_{loc}(\R^n)$.)  We will also
assume that
$V\ge 0$ (the case when $V$ is semi-bounded below by another constant is easily reduced to
the case when $V\ge 0$ for the discreteness of spectrum results). Then we can define $H_{a,V}$ 
as the operator defined by
the closure of this quadratic form. This closure is well defined \cite{Leinfelder-Simader}. 

We will say that  $H_{a,V}$
has a \textit{discrete spectrum} if 
its spectrum consists of isolated eigenvalues of finite
multiplicities. 
It follows that the only accumulation point of
these eigenvalues can be $+\infty$. Equivalently we may say that $H_{a,V}$  
has a compact resolvent.

Our first goal is to provide necessary and sufficient conditions  for the discreteness of the
spectrum of $H_{a,V}$. We will write  $\sigma=\sigma_d$ instead of the statement that the 
spectrum of $H_{a,V}$ is discrete. 

Let us recall some facts concerning the Schr\"odinger
operator $H_{0,V}=-\De+V$ without magnetic field (i.e.\ the operator
(\ref{E:Schroed}) with $a=0$).

It is a  classical result of K.~Friedrichs \cite{Friedrichs}
(see also  e.g.\ \cite{Reed-Simon}, Theorem XIII.67, or 
\cite{Berezin-Shubin}, Theorem 3.1) that
the condition
\begin{equation*}\label{E:V-to-infty}
V(x)\to+\infty\quad\hbox{as}\quad x\to\infty
\end{equation*}
implies  $\sigma=\sigma_d$ (for $H_{0,V}$).

A.~Molchanov \cite{Molchanov} found  a necessary and sufficient condition for
the discreteness of spectrum. 
It is formulated in terms of the Wiener  capacity. The capacity of a compact set $F$  will be 
denoted $\capa(F)$ (see Section \ref{S:capacity} for the definition and 
\cite{Edmunds-Evans, Kondrat'ev-Shubin, Mazya8} for  necessary properties of the capacity, 
expositions of Molchanov's work and more general results).  

Let $B(x,r)$ denote the open ball in $\R^n$ with the radius $r>0$ and the center at $x$,
$\bB(x,r)$ denote the corresponding closed ball.  

In case $n=2$ the capacity of a set $F\subset \bB(x,r)$  is always taken relative to
a ball $B(x,2r)$. The value of $r$ is usually clear from the context. 
In case $n\ge 3$ such a definition would be
equivalent to the usual Wiener capacity  (relative to $\R^n$). 

In case $n=2$ we can also use capacities of sets $F\subset\bB(x,r)$
with respect to the ball  $B(x,R)$ where $r\in (0,R/2)$ and $R>0$ is fixed, but this 
complicates some formulations.  

Similarly we can use closed cubes (squares if $n=2$)
$Q_d$, where $d>0$ means the length of the edge and the edges are assumed to be parallel to the
coordinate axes. The interior of $Q_d$ will be denoted $\overset{\circ}Q_d$.
In this paper we prefer to use cubes instead of balls, but balls are more
convenient in case of manifolds. 
In case $n=2$ the capacity of a compact set $F\subset Q_d$ will be always defined relative to
$\overset{\circ}Q_{2d}$, where $Q_d$ and $Q_{2d}$ have the same center.

Let us define the Molchanov functional
\begin{equation}\label{E:Mc}
M_c(Q_d;V)=
\inf_F\left\{\left.\int_{Q_d\setminus F} V(x)dx\right|\capa (F)\le c\,\capa (Q_d)\right\}.
\end{equation}
Here we will always assume that $0<c<1$. Due to the standard properties of the capacity, the infimum in 
\eqref{E:Mc} will not change if we only restrict it to the sets $F$ which are 
closures of  open subsets of $Q_d$ with a smooth boundary.

A.~Molchanov proved that there exists $c=c_n>0$ such that $H_{0,V}$ has a discrete spectrum
if and only if for every $d>0$
\begin{equation}\label{E:M-cond}
M_c(Q_d;V)\to+\infty \quad \hbox{as}\quad Q_d\to\infty,
\tag{$M_c$}
\end{equation}
where $Q_d\to\infty$ means that the center of the cube $Q_d$ 
goes to infinity (with $d$ fixed). He actually established
this result with a specific constant $c_n$ (see also \cite{Kondrat'ev-Shubin}), namely, 
$c_n=(4n)^{-4n}(\capa(Q_1))^{-1}$ for $n\ge 3$, but it is by no means precise and we will 
not be interested in the precise value of this constant (it seems beyond the reach
of the existing technique). 

The case $n=2$ was not discussed in
\cite{Molchanov}, though it can be  covered by the same methods with minor modifications.

Note that $(M_c)$ implies $(M_{c'})$ for every $c'<c$. 
The arguments in \cite{Molchanov} actually show that it suffices to assume that
$(M_c)$ is satisfied for all sufficiently small $c>0$.
Hence we can equivalently formulate a necessary and sufficient condition 
of the discreteness of spectrum for $H_{0,V}$ by writing that
$(M_c)$ is satisfied for all
$c\in(0,c_0)$ with a positive $c_0$. 

Note also that $\capa(\bB(x,r))$ can be explicitly calculated.
It equals $c_nr^{n-2}$ 
(with a different $c_n>0$).
The capacity of a cube $Q_d$ is $c_nd^{n-2}$ 
(with yet another $c_n>0$).
Hence in the formulation
of the Molchanov condition $(M_c)$ we can replace $\capa(Q_d)$
by $d^{n-2}$. 

A simple argument given in \cite{Avron-Herbst-Simon} (see also
Corollary 1.4 in \cite{Kondratiev-Shubin}) shows that if
$H_{0,V}$ has a discrete
spectrum, then the same is true for $H_{a,V}$ whatever the
vector potential $a$. Therefore the condition (\ref{E:M-cond})
together with $V\ge 0$ is sufficient for the
discreteness of spectrum of $H_{a,V}$. This means that a magnetic
field can only improve the situation from our point of view.
Papers by J.~Avron, I.~Herbst and B.~Simon
\cite{Avron-Herbst-Simon}, Y.~Colin de Verdi\`ere \cite{Colin-de-Verdiere}, 
A.~Dufresnoy \cite{Dufresnoy}
and A.~Iwatsuka \cite{Iwatsuka1} provide some quantitative
results which show that even in case $V=0$ the magnetic field 
can make the spectrum discrete. (This situation is called
\textit{magnetic bottle}.) 

The results of 
\cite{Avron-Herbst-Simon, Dufresnoy, Iwatsuka1},  
were improved in \cite{Kondratiev-Shubin}.  
In particular,
some sufficient conditions for the spectrum of $H_{a,V}$
to be discrete were given. 
The capacity was added into the picture, so in most cases these conditions
become necessary and sufficient in case when there is no magnetic
field, i.e.\ when $a=0$. Also both electric and magnetic fields
were made to work together to achieve
the discreteness of spectrum. 

However no  necessary and sufficient conditions of the discreteness 
of the spectrum
with both fields present were provided in \cite{Kondratiev-Shubin}. 
Here we will give such conditions which
actually separate the influence of the electric and magnetic fields.
If the magnetic field is absent then our conditions turn into the Molchanov condition 
\eqref{E:M-cond} or into some weaker conditions,  improving Molchanov's 
sufficiency result.  

We will need the bottoms 
$\la(G;H_{a,V})$ and $\mu(G;H_{a,V})$ of 
Dirichlet and Neumann spectra for
the operator $H_{a,V}$ in an open set $G\subset\R^n$.
They are defined in terms of its quadratic form $h_{a,V}$
as follows (see e.g.\ \cite{Courant-Hilbert}, \cite{Kato1}):

\begin{equation}\label{E:lambda}
\la(G;H_{a,V})=\underset{u}{\inf} \left\{\frac{h_{a,V}(u,u)_G}
{(u,u)_G},\;
u\in C_c^\infty(G)\setminus \{0\}\right\}\;,
\end{equation}

\begin{equation}\label{E:mu}
\mu(G;H_{a,V})=\underset{u}{\inf} \left\{\frac{h_{a,V}(u,u)_G}
{(u,u)_G},\;
u\in (C^\infty(G)\setminus\{0\})\cap L^2(G)
\right\}\;,
\end{equation}
where in both cases $h_{a,V}(u,u)_G$ is given by the formula
\eqref{E:haV} with the integrals over $G$ (instead of $\R^n$) i.e.
\begin{equation*}
h_{a,V}(u,u)_G=\int_{G}(|\nabla_a u|^2+V|u|^2)dx,
\end{equation*}
and $(u,u)_G$ means square of the $L^2$-norm of $u$ in $G$.
However in the future we will often skip the subscript $G$
since it will be clear from the context which $G$ is used.

We will also use these notations for $G=Q_d$ in which case $\la(Q_d;H_{a,V})$ is understood
as $\la(\overset{\circ}Q_d;H_{a,V})$, whereas $\mu(Q_d;H_{a,V})$ can be understood as 
$\mu(\overset{\circ}Q_d;H_{a,V})$ as well as directly by the formula \eqref{E:mu}
(i.e. with the use of functions $u$ which are $C^\infty$ on the closed cube) which gives 
the same result.

In both \eqref{E:lambda} and \eqref{E:mu} 
we can also use locally Lipschitz test functions instead of $C^\infty$ functions $u$,
which does not change the result. (Of course we should take functions with compact support 
in $G$ in case of $\la(G;H_{a,V})$.)

We will also need the quantity
\begin{equation}\label{E:mu0}
\mu_0=\mu_0(Q_d)=\mu_0(Q_d;a)=\mu(Q_d;H_{a,0}),
\end{equation}
which we will call the {\it local energy of the magnetic field} (in $Q_d$). 
Here the first three terms are defined by the last one, but we will use 
the shorter notations when the choice of $Q_d$ and $a$ is clear from the context.
Obviously $\mu_0\ge 0$. Also, $\mu_0$ is {\it gauge invariant} i.e. 
\begin{equation*}\label{E:gauge}
\mu_0(Q_d;a)=\mu_0(Q_d;a+d\phi),
\end{equation*}
as soon as $a, a+d\phi\in L^\infty_{loc}(Q_d)$, $\phi$ is a locally Lipschitz function, 
and $a$ is identified with the 1-form
\begin{equation*}\label{E:a-form}
a=\sum_{j=1}^n  a_j dx^j.
\end{equation*}
Therefore $\mu_0(Q_d;a)$ depends only on the magnetic field $B=da$
which is understood as a 
2-form with distributional coefficients.
It is easy to see that $\mu_0(Q_d;a)$ vanishes if and only if 
$B$ vanishes on $\overset{\circ}Q_d$.  This justifies calling
$\mu_0$ local energy of the magnetic field.

We will also use a {\it normalized local energy of the magnetic field} in $Q_d$ 
defined as
\begin{equation}\label{E:tmu0}
\tmu_0=\tmu_0(Q_d)=\tmu_0(Q_d;a)=\mu_0 d^2.
\end{equation}

\begin{definition}\label{D:F}
{\rm
A class ${\cF}$ consists of functions $f:[0,+\infty)\to (0,+\infty)$
which are continuous and decreasing on $[0,+\infty)$.

A class $\cG$ consists of functions $g:(0,d_0)\to(0,+\infty)$ such that $g(\tau)\to 0$
as $\tau\to 0$ and $(g(d))^{-1}d^2\le 1$  for all $d\in (0,d_0)$.

The pair $(f,g)\in \cF\times\cG$ is called {\it $n$-admissible} if $f$ satisfies
the inequality
$f(t)\le f_n(t)$ for all $t\ge 0$, where 
\begin{equation}\label{E:f_n}
f_n(t)=(1+t)^{(2-n)/2} \ \ {\rm if}\  \ n\ge 3,\qquad f_2(t)=(1+\log(1+t))^{-1}.
\end{equation}
}
\end{definition}

Now we can formulate our main result about the discreteness of spectrum.

\begin{theorem}\label{T:discr} Let us assume that $a\in L^\infty_{loc}(\R^n)$. 
There exists
$c_n>0$ such that for every $n$-admissible pair $(f,g)$ the following conditions on $H_{a,V}$ are equivalent:

\ms
$({\rm a})$
The spectrum of $H_{a,V}$ is discrete.

\ms

$({\rm b}_{f,g})$ 
There exists $d_0>0$ such that for every $d\in (0,d_0)$
\begin{equation}\label{E:b_f}
\mu_0(Q_d)+d^{-n}M_\ga(Q_d;V)\to +\infty
 \quad {\rm as}\quad Q_d\to\infty,
\end{equation}
where 
\begin{equation}\label{E:gamma}
\gamma=\gamma(\mu_0,d)=c_nf(\tmu_0)g(d)^{-1}d^2. 
\end{equation}

\ms
$({\rm c}_{f,g})$ 
There exists $d_0>0$ such that for every $d\in(0,d_0)$
\begin{equation}\label{E:cfg}
\liminf_{Q_d\to\infty} \left(\mu_0(Q_d)+d^{-n}M_\gamma(Q_d;V)\right)\ge g(d)^{-1},
\end{equation}
where $\gamma$ is as in \eqref{E:gamma}.
\end{theorem}

Note that   
$f(\tmu_0)=f(\mu_0 d^2)$ is decreasing in $\mu_0$ and tends to $0$ as $\mu_0\to\infty$
(with $d$ fixed).
So the condition on $V$ is weaker at the places where the local energy of the 
magnetic field is larger. 

\begin{remark}\label{R:no-magn}{\rm
Assuming that  the magnetic field is absent ($a=0$, $H_{a,V}=H_{0,V}=-\Delta +V$) we obtain  
$c_nf(\tmu_0)=c_n f(0)=c>0$. Now taking $g(d)=d^2$ we see that the condition 
\eqref{E:b_f} becomes  the Molchanov condition $(M_c)$.
So Theorem \ref{T:discr} strengthens Molchanov's theorem \cite{Molchanov}
which claims the equivalence of $({\rm a})$ and $({\rm b}_{f,g})$ for this particular case.
}
\end{remark}

\begin{corollary}\label{C:f-equiv}
All conditions $({\rm b}_{f,g})$, $({\rm c}_{f,g})$, taken for different
$n$-admissible pairs $(f,g)$ are equivalent.
\end{corollary}

In particular, this Corollary applied in case $a=0$ (no magnetic field) gives
an equivalence of different conditions on the scalar potential $V\ge 0$. This seems to be 
a new purely function-theoretic property of capacity.

\ms
The following corollaries provide examples of more explicit  necessary and separately 
sufficient conditions which easily follow from Theorem \ref{T:discr}.

\begin{corollary}\label{C:necessary}
Let us assume that the spectrum of $H_{a,V}$ is discrete. Then
for every fixed $d>0$ 
\begin{equation}\label{E:cor-nec}
\mu_0(Q_d)+\frac{1}{d^n}\int_{Q_d}V(x)dx\to +\infty
 \quad as \quad Q_d \to \infty.
\end{equation}
\end{corollary}

The condition \eqref{E:cor-nec} corresponds to the case $\gamma\equiv 0$ in $({\rm b}_{f,g})$ 
in Theorem \ref{T:discr}. It is known that it is not sufficient for the discreteness of the spectrum,
even in the case when there is no magnetic field \cite{Molchanov}.

\begin{corollary}\label{C:sufficient}
Let us assume that there exist $c>0, d_1>0$
such that for every fixed $d\in (0,d_1)$ 
\begin{equation}\label{E:cor-suf}
\mu_0(Q_d)+d^{-n}M_c(Q_d;V)\to +\infty
 \quad as \quad Q_d \to \infty.
\end{equation}
Then the spectrum of $H_{a,V}$ is discrete.
\end{corollary}

It follows from Theorem \ref{T:precision} below that the condition
\eqref{E:cor-suf} is not necessary for the discreteness of spectrum 
of $H_{a,V}$.

Sufficient conditions (for $\sigma=\sigma_d$)  which do not include capacity,
can be obtained if the capacity is replaced by the Lebesgue measure 
in the restriction on $F$  in the definition of $M_c(Q_d;V)$ 
-- see Section 6.1 in \cite{Kondrat'ev-Shubin} for a more
detailed argument.

Other, more effective sufficient conditions (which do not include $\mu_0$)
and related results (in particular, asymptotics of eigenvalues under appropriate conditions) 
can be found in 
\cite{Colin-de-Verdiere, Dufresnoy, Fefferman, Helffer-Mohamed, Helffer-Nourrigat-Wang, Ivrii,
Iwatsuka1, Iwatsuka2, Kondratiev-Shubin, Levendorskii, Mohamed-Raikov, Shigekawa, Tamura}.

Some necessary and sufficient conditions of discreteness of spectrum for 
the Schr\"odinger operators can be obtained by considering them as 
1-dimensional Schr\"odinger operators with operator coefficients
(see e.g. \cite{Maslov, Bruning} and references in \cite{Bruning}). An interesting feature 
of this approach is that it allows to consider operators whose potentials are not necessarily 
semi-bounded below.

\ms
The following Theorem shows that the conditions on $f$ in Theorem \ref{T:discr}
are almost precise.

\begin{theorem}\label{T:precision}
There exists an operator $H_{a,V}$ with a discrete spectrum and with the following property.
Let $f:[0,+\infty)\to(0,1)$ be a decreasing function, such that in case $n\ge 3$
\begin{equation}\label{E:not-adm-ge-3}
f(t)=(1+t)^{\frac{2-n}{2}}h(t), 
\end{equation}
and in case $n=2$
\begin{equation}\label{E:not-adm-2}
f(t)=(1+\log(1+t))^{-1}h(t), 
\end{equation}
where in both cases $h(t)\to +\infty$ as $t\to +\infty$.
Then, for every fixed $d>0$,
the condition \eqref{E:b_f} with $\ga=f(\mu_0 d^2)$
is not satisfied. 
So the condition \eqref{E:b_f} with the function $f$ having the form given above,
is not necessary for the discreteness of spectrum, whatever $g$ and $c_n$. 
In particular, the exponents in 
\eqref{E:f_n} are the best possible.
\end{theorem}

\ms 
Now we will give a positivity criterion for the operators $H_{a,V}$. We will say that 
such an operator is {\it strictly positive} if $H_{a,V}\ge\eps I$ for some $\eps>0$,
or, equivalently, that its spectrum is in $[\eps,\infty)$ for some $\eps>0$. If $V\ge 0$,
then this is equivalent to saying that $0$ is not in the spectrum of $H_{a,V}$. 

\begin{theorem}\label{T:pos} 
Let us assume that $V\ge 0$. There exist positive constants $c_n,\tc_n$  such that 
the following conditions on $H_{a,V}$ are equivalent: 
 
$({\rm a})$ \quad $H_{a,V}$ is strictly positive.

$({\rm b})$ There exist positive constants $c,d_1,d$ 
such that for  every cube $Q_d\subset \R^n$  
\begin{equation}\label{E:M-pos}
\mu_0(Q_d)+d^{-n}M_c(Q_d;V)\ge \frac{1}{d_1^2}.
\end{equation} 

$({\rm c})$ There exist positive constants $d_1,d$   
such that for  every cube $Q_d\subset \R^n$  
\begin{equation}\label{E:M-pos-n}
\mu_0(Q_d)+d^{-n}M_{c_n}(Q_d;V)\ge \frac{1}{d_1^2}.
\end{equation}

$({\rm d})$ There exist positive constants $c, \tc, d_2$
such that for  every $d>d_2$  and every cube $Q_d\subset \R^n$  
\begin{equation}\label{E:M-pos-d}
\mu_0(Q_d)+d^{-n}M_c(Q_d;V)\ge \frac{\tc}{d^2}.
\end{equation} 

$({\rm e})$ There exist $d_2>0$ 
such that for  every $d>d_2$  and every cube $Q_d\subset \R^n$  
\begin{equation}\label{E:M-pos-dn}
\mu_0(Q_d)+d^{-n}M_{c_n}(Q_d;V)\ge \frac{\tc_n}{d^2}.
\end{equation}     

\end{theorem}

In case when there is no magnetic filed (i.e. $a=0$, $H_{a,V}=H_{0,V}=-\Delta+V$) this 
theorem is essentially contained in \cite{Mazya8}, Sect. 12.5.

\ms
\begin{remark}\label{R:gauge}
{\rm 
The discreteness of spectrum and strict positivity are gauge invariant. More precisely, if we replace  
$a\in L^\infty_{loc}(\R^n)$ by another magnetic potential $a'\in L^\infty_{loc}(\R^n)$ which has the form
$a'=a+d\phi$, then the spectrum does not change, i.e. the spectra of $H_{a,V}$ and $H_{a',V}$
coincide (see \cite{Leinfelder}).  (Here $\phi$ is a locally Lipschitz function.)  
So in fact the spectrum depends not on
the magnetic potential $a$ itself but on the magnetic field $B=da$. 
}
\end{remark}

\begin{remark}\label{R:manifolds}
{\rm
Theorem \ref{T:discr} holds on every manifold of bounded geometry, with cubes 
replaced by balls in the formulation
(see \cite{Kondrat'ev-Shubin}  and Section 6 in \cite{Kondratiev-Shubin} for necessary adjustments 
which should be done to treat the more general case compared with the case of operators on $\R^n$). 
However it is not at
all clear how to extend Theorem
\ref{T:pos} to this case. 
}
\end{remark}

\begin{remark}
{\rm 
In Section \ref{S:domains} we will formulate results which extend Theorems \ref{T:discr}
and \ref{T:pos} and their Corollaries to the case when the operator $H_{a,V}$
is considered in $L^2(\Om)$ for an arbitrary open set $\Om\subset \R^n$ with the Dirichlet
boundary conditions on $\pa\Om$. Note that the discreteness of spectrum and strict positivity 
in this case may be influenced or even completely determined by the geometry of $\Om$. In particular, 
the results are non-trivial even for the pure Laplacian $H_{0,0}=-\De$.  
}
\end{remark}

\section{Preliminaries}\label{S:capacity}

In this section we will list some important technical tools which will be used later.
They were actually useful even in case of vanishing magnetic field (see \cite{Mazya8}), when they
provide simpler proofs and stronger versions for the Molchanov discreteness of spectrum
criterion, as well as for the Maz'ya strict positivity criterion for usual
Schr\"odinger operators with non-negative scalar potentials.

For every subset $\Omega\subset\R^n$ denote by $\Lip(\Omega)$ the space of 
(complex-valued) functions satisfying
the uniform Lipschitz condition in $\Omega$, and by $\Lip_c(\Omega)$ the subspace in $\Lip(\Omega)$
of all functions with compact support in $\Omega$ (this will be only used when $\Omega$ is open).
By $\Lip_{loc}(\Omega)$ we will denote the set of functions on (an open set) $\Omega$ which are
Lipschitz on any compact subset $K\subset\Omega$.

If $F$ is a compact subset in  an open set $\Omega\subset{\R}^n$, then the 
Wiener capacity of $F$ relatively to $\Omega$
is defined as
\begin{equation}\label{E:cap-def}
{\capa}_\Omega (F)=\inf\left\{\left.
\int_{\R^n}|\nabla u(x)|^2 dx\,\right|\; u\in\Lip_c(\Om), u|_F=1\right\}.
\end{equation}

We will also use the notation $\capa(F)$ for $\capa_{\R^n}(F)$
if $F\subset \R^n$, $n\ge 3$, and for 
$\capa_{\overset{\circ}Q_{2d}}(F)$ if $F\subset Q_d\subset \R^2$,  
where the squares $Q_{d}$ and $Q_{2d}$ have the same center and
the edges parallel to  the coordinate axes in $\R^2$.

Note that if we allow only real-valued functions $u$ in \eqref{E:cap-def}, then
the infimum will not change. To see this it suffices to note that
$|\nabla |u||\le |\nabla u|$ a.e. (almost everywhere) for every complex-valued Lipschitz function.
Moreover, the infimum does not change if we  restrict ourselves to the Lipschitz 
functions $u$ such that $0\le u\le 1$  everywhere (see e.g. \cite{Mazya8}, Sect. 2.2.1).

The following  Lemmas are particular cases of much more general results from \cite{Mazya8}. 
We supply the simplified formulations for the convenience  of the readers. 

\begin{lemma}\label{L:Mazya1} {\rm (\cite{Mazya8}, Theorem 10.1.2, part 1).}
There exists $C_n>0$ such that the following inequality holds for every complex-valued
function $u\in \Lip(Q_d)$  which
vanishes on a compact set $F\subset Q_d$ (but is not identically zero on $Q_d$):
\begin{equation}\label{E:cap-above}
\capa(F)\le \frac{C_n\int_{Q_d} |\nabla u(x)|^2dx}{d^{-n}\int_{Q_d} |u(x)|^2 dx}\;.
\end{equation}
\end{lemma}

\begin{lemma}\label{L:Mazya2} \rm{(\cite{Mazya8}, Lemma 12.1.1).}
 Let $V\in L^1_{loc}(\R^n)$, $V\ge 0$.
For every $u\in \Lip(Q_d)$ and  $\gamma>0$
\begin{equation}\label{E:Mazya(2)}
\int_{Q_d}|u|^2 dx\le {C_nd^2\over \gamma}  \int_{Q_d}|\nabla u|^2 dx
+{4d^n\over M_\ga(Q_d;V)}\
\int_{Q_d}V|u|^2 dx,
\end{equation}
(The last term is declared to be $+\infty$ if its denominator vanishes.)
\end{lemma}

\begin{remark}\label{R:main-ineq}
{\rm 
Both Lemmas \ref{L:Mazya1} and \ref{L:Mazya2} hold also if we replace
$\nabla$ by $\nabla_a$. Indeed, we can first apply the inequalities \eqref{E:cap-above}
and \eqref{E:Mazya(2)} to $|u|$ and then use the diamagnetic inequality $|\nabla|u||\le |\nabla_au|$
(see e.g.\ \cite{Kato2, Lieb-Loss, Simon76}).
}
\end{remark}

\bigskip
The following lemma is somewhat inverse to Lemma \ref{L:Mazya1}.
It follows from part 2 of Theorem 10.1.2 in \cite{Mazya8}.

\begin{lemma}\label{L:cap-F-below} There exists positive $c_n, c_n', c_n''$ such that 
for every compact subset $F'\subset Q_d$ satisfying 
\begin{equation}\label{E:cap-F'-cn}
\capa(F')\le c_n\capa(Q_d),
\end{equation}
there exists $\psi\in\Lip(Q_d)$ with the following properties:
$0\le\psi\le 1$,   $\psi=0$ in a neighborhood of $F'$,
\begin{equation}\label{E:cap-F'-psi-below}
\capa(F')\ge c'_n\int_{Q_{d}}|\nabla \psi|^2dx
\end{equation}
and
\begin{equation}\label{E:psi2-below}
d^{-n}\int_{Q_d}\psi^2 dx\ge \frac{1}{4},
\end{equation}
hence
\begin{equation}\label{E:cap-F'-below}
\capa(F')\ge \frac{c''_n\int_{Q_d}|\nabla\psi|^2dx}{d^{-n}\int_{Q_d}\psi^2dx}\;.
\end{equation}
\end{lemma}

\ms
For the convenience of the reader we provide self-contained  proofs of the lemmas above in 
Appendix to this paper.

\section{Discreteness of spectrum: sufficiency.}\label{S:sufficient}

In this section we will consider operators $H_{a,V}$ with $V\in L^1_{loc}(R^n)$, 
$V\ge 0$ and $a\in L^\infty_{loc}(\R^n)$.

We will start with the following proposition which  gives 
a general (albeit complicated) sufficient condition for the discreteness of spectrum. 

\begin{proposition}\label{P:general-suff}
Given an operator $H_{a,V}$, let us assume that the following condition is satisfied:
\begin{align*}
&\exists\; \eps_0>0,\ \forall\;\eps\in (0,\eps_0),\ \exists\; d=d(\eps)>0,\; R=R(\eps)>0, 
\forall\; Q_d \ \text{with} \\ 
&Q_d\cap (\R^n\setminus B(0,R))\ne\emptyset,\  
\exists \ga=\ga(\mu_0, d, \eps)\ge 0, \ \text{such that}
\end{align*}
\begin{equation}\label{E:2-main-ineq}
\mu_0+\frac{\ga}{C_n d^2}\ge \eps^{-1}\quad \text{and}\quad \mu_0+d^{-n}M_\ga(Q_d;V)\ge \eps^{-1},
\end{equation}
where $\mu_0=\mu_0(Q_d)$, $C_n$ is the constant from \eqref{E:Mazya(2)}. Then $\sigma=\sigma_d$.
\end{proposition}

{\bf Proof.} We can assume without loss of generality that $V\ge 1$.
Define
\begin{equation}\label{E:cL}
\cL=\left\{u\left|u\in C_c^\infty(\R^n), 
\int_{\R^n}(|\nabla_a u|^2+V|u|^2)dx\le 1\right.\right\}.
\end{equation}
By the standard functional analysis argument (see e.g. Lemma 2.3 in \cite{Kondrat'ev-Shubin}) 
the spectrum  of $H_{a,V}$ is discrete if and only if $\cL$ is precompact in $L^2(\R^n)$,
which in turn holds if and only if $\cL$ has ``small tails", 
i.e. for every $\eps>0$ there exists $R>0$ such that 
\begin{equation}\label{E:smalltails}
\int_{\R^n\setminus B(0,R)}|u|^2 dx\le \eps \quad \text{for all} \quad u\in\cL.
\end{equation}
This will hold if we establish that there exists $d>0$ such that 
\begin{equation}\label{E:cube-tail}
\int_{Q_d}|u|^2 dx\le \eps \int_{Q_d}(|\nabla_a u|^2+V|u|^2)dx, 
\end{equation}
for all cubes $Q_d$ such that $Q_d\cap(\R^n\setminus B(0,R))\ne\emptyset$.

To prove \eqref{E:cube-tail} note first that if $\gamma=0$ then $\mu_0\ge \eps^{-1}$
due to the first inequality in \eqref{E:2-main-ineq}, hence
\eqref{E:cube-tail} follows from the definition of $\mu_0$ (even if we skip the term
with $V$ in the right-hand side). So from now we will assume that $\gamma>0$.

Let us look at the inequality
\begin{equation}\label{E:Mazya(2)-magn}
\int_{Q_d}|u|^2 dx\le {C_nd^2\over \gamma}\int_{Q_d}|\nabla_a u|^2 dx
+{4d^n\over M_\ga(Q_d;V)}\
\int_{Q_d}|u|^2 Vdx
\end{equation}
(see Lemma \ref{L:Mazya2} and Remark \ref{R:main-ineq}).
For every fixed $\eps>0$ we can divide  all cubes $Q_d$ 
into the following two types: 

\ms
Type I:\qquad $\mu_0(Q_d)>(2\eps)^{-1}$;

\ms
Type II:\qquad  $\mu_0(Q_d) \le (2\eps)^{-1}$.

\ms
For a Type I cube $Q_d$ the inequality \eqref{E:cube-tail} holds with $2\eps$
instead of $\eps$, as was explained above.

For a Type II cube it follows from the conditions \eqref{E:2-main-ineq} that
\begin{equation*}
{C_nd^2\over \gamma}\le 2\eps, \qquad {4d^n\over M_\ga(Q_d;V)}\le 8\eps,
\end{equation*}
so the inequality \eqref{E:cube-tail} follows with $8\eps$ instead of $\eps$. $\square$

Instead of requiring that the conditions of Proposition \ref{P:general-suff} satisfied 
for all $\eps\in(0,\eps_0)$, it suffices to require it for a sequence $\eps_k\to +0$.
Keeping this in mind we can replace the dependence $d=d(\eps)$ by the inverse
dependence $\eps=g(d)$, so that $g(d)>0$ and $g(d)\to 0$ as $d\to +0$ (and here we can also
restrict to a sequence $d_k\to +0$). This leads to the following
\begin{proposition}\label{P:general-suff-d}
Given an operator $H_{a,V}$ with $V\ge 0$, let us assume that the following condition is satisfied:
\begin{align*}
&\exists\; d_0>0,\ \forall\; d\in (0,d_0),\ \exists\; R=R(d)>0, 
\forall\; Q_d \ \text{with} \\ 
&Q_d\cap (\R^n\setminus B(0,R))\ne\emptyset,\  
\exists \ga=\ga(\mu_0, d)\ge 0, \ \text{such that}
\end{align*}
\begin{equation}\label{E:2-main-ineq-d}
\mu_0+\frac{\ga}{C_n d^2}\ge g(d)^{-1}\quad \text{and}\quad \mu_0+d^{-n}M_\ga(Q_d;V)\ge g(d)^{-1},
\end{equation}
where $\mu_0=\mu_0(Q_d)$, $C_n$ is the constant from \eqref{E:Mazya(2)},
$g(d)>0$ and $g(d)\to 0$ as $d\to +0$. Then $\sigma=\sigma_d$.
\end{proposition}

\begin{proposition}\label{P:sufficient} Let us assume that $V\ge 0$,
$f\in\cF$, $g\in\cG$ (in the notations of Definition \ref{D:F}) and
one of the conditions $({\rm b}_{f,g})$,
$({\rm c}_{f,g})$ from Theorem \ref{T:discr} is satisfied.
Then the spectrum of $H_{a,V}$ is discrete.
\end{proposition}

{\bf Proof.} Clearly, 
$({\rm b}_{f,g})$ implies
$({\rm c}_{f,g})$.
So it remains to prove that $({\rm c}_{f,g})$  implies that $\sigma=\sigma_d$.
To this end it is sufficient to prove that it implies that the conditions 
of Proposition \ref{P:general-suff-d} are satisfied. 

Note that it suffices to establish that the inequalities \eqref{E:2-main-ineq-d} hold 
with an additional positive constant factor, independent on $d$ (but possibly dependent on
$f,g$), in the right hand sides. 

Clearly, the second inequality in \eqref{E:2-main-ineq-d}, with an additional factor $1/2$
in the right hand side, is satisfied 
for distant cubes $Q_d$ due to \eqref{E:cfg}. So we need only to take care for the first
inequality in \eqref{E:2-main-ineq-d}. It obviously holds if $\mu_0\ge g(d)^{-1}$.

On the other hand, if we assume that $\mu_0 \le g(d)^{-1}$, then
\begin{equation*}
f(\mu_0 d^2)\ge f(g(d)^{-1}d^2),
\end{equation*}
hence
\begin{equation*}
\frac{\gamma}{C_n d^2}=\frac{c_n}{C_n} f(\mu_0 d^2)g(d)^{-1}\ge
\frac{c_n}{C_n}f(g(d)^{-1}d^2)g(d)^{-1} \ge \frac{c_n}{C_n}f(1)g(d)^{-1},
\end{equation*}
because $g(d)^{-1}d^2\le 1$ according to Definition \ref{D:F}. Therefore we can apply Proposition 
\ref{P:general-suff-d}. 
$\square$

\begin{remark}\label{R:free-f}
{\rm
No domination requirement (like $f\le f_n$ in Definition \ref{D:F}) is imposed on $f$
in Proposition \ref{P:sufficient}.
}
\end{remark}

\begin{remark}\label{R:tiling}
{\rm
It is clear from the proof that to establish the discreteness
of spectrum of an operator $H_{a,V}$, it suffices to check
the condition $({\rm b}_{f,g})$ (or $({\rm c}_{f,g})$)
from Theorem \ref{T:discr} for every $d\in (0,d_0)$
on the cubes $Q_d$ which form a tiling of $\R^n$
(instead of all cubes $Q_d$).
}
\end{remark}

\begin{remark}\label{R:M-improve}
{\rm
Let us consider the case of vanishing magnetic field ($a\equiv 0$) and take $g(d)=d^s$ with $0<s <2$.
Then the conditions $({\rm b}_{f,g})$, $({\rm c}_{f,g})$  provide sufficient conditions 
for the discreteness of spectrum of the Schr\"odinger operator $H_{0,V}=-\Delta+V$ which are much better
than the Molchanov condition $(M_c)$ which corresponds to the condition $({\rm b}_{f,g})$ with $g(d)=d^2$.
The conditions $({\rm b}_{f,d^s})$ in this case impose weaker requirements on the capacity
of negligible sets for small $d$. With the same requirements on the negligible sets
the condition $({\rm c}_{f,d^s})$ goes even further: it does not require the functional $M_\gamma(Q_d;V)$
to go to infinity for fixed $d$, it only requires it to become large for distant cubes and small $d$. 
}
\end{remark}

\section{Discreteness of spectrum: necessity.}\label{S:nec}

We will use the notations from Section \ref{S:Intro}. We impose here the same restrictions on 
$H_{a,V}$ as in Section \ref{S:sufficient}, i.e. $V\in L^1_{loc}(\R^n)$, 
$V\ge 0$, $a\in L^\infty_{loc}(\R^n)$.
Let us fix an
arbitrary
$d_0>0$. We need to prove that the discreteness of spectrum for $H_{a,V}$
implies the condition $({\rm b}_{f,g})$ in Theorem \ref{T:discr}.
This will follow from 

\begin{proposition}\label{P:nec}
 There exist $c=c_n>0$, $C=C_n>0$  
such that for 
every operator $H_{a,V}$ with $V\ge 0$
and every cube $Q_d$
\begin{equation}\label{E:nec-ineq}
\mu(Q_d;H_{a,V})\le C E\left(1+\frac{1}{f_n(\tmu_0)d^{n-2}}M_{c f_n(\tmu_0)}(Q_d;V)\right),
\end{equation}
where $E=\mu_0(Q_d)+d^{-2}$, $\tmu_0$ is defined by \eqref{E:tmu0}, 
and $f_n$ is defined by \eqref{E:f_n}.
\end{proposition}

\bigskip
{\bf Proof of Theorem \ref{T:discr}.} Clearly $({\rm b}_{f,g})$ implies $({\rm c}_{f,g})$.
The sufficiency of the condition $({\rm c}_{f,g})$  for the discreteness of spectrum was proved  in Section
\ref{S:sufficient}. So we only need to prove that $\sigma=\sigma_d$ implies $({\rm b}_{f,g})$
for every $n$-admissible pair $f,g$ (see Definition \ref{D:F}). It is sufficient to consider the special case $f=f_n$,
$g(d)=d^2$ because this case corresponds to the maximal allowed value of $\gamma(\mu_0,d)$, therefore to the strongest
possible condition $({\rm b}_{f,g})$ among all possible $n$-admissible pairs $(f,g)$.

So let us assume that
$H_{a,V}$ has a discrete spectrum. We need to prove that the condition  
$({\rm b}_{f,g})$ holds for $f=f_n$, $g(d)=d^2$. For brevity sake denote this condition by $(N)$.

According to the  Localization Theorem 1.2 in \cite{Kondratiev-Shubin} 
it follows from the discreteness of spectrum that
\begin{equation}\label{E:mu-to-infty}
\mu(Q_d;H_{a,V})\to +\infty\quad\text{as}\quad Q_d\to\infty,
\end{equation}
for every fixed $d>0$. 
 This implies that the right hand side of \eqref{E:nec-ineq} tends to $+\infty$
as $Q_d\to\infty$ with any fixed $d>0$. This implies that the condition 
$(N)$ is satisfied. Indeed, if $(N)$ does not hold for some $d>0$,
then 
there exists 
a sequence of cubes $Q_d\to\infty$ such that 
\begin{equation*}
E+d^{-n}M_{cf_n(\tmu_0)}(Q_d;V)\le C
\end{equation*}
along this sequence. But then both terms in the left hand side are bounded,
hence the right hand side of \eqref{E:nec-ineq} is bounded,
which contradicts \eqref{E:mu-to-infty}. $\square$

\bigskip
Now we will start our proof of Proposition \ref{P:nec}. Let us choose $u\in\Lip(Q_d)$
such that
\begin{equation}\label{E:ha0-ineq}
h_{a,0}(u,u)=\int_{Q_d}|\nabla_a u|^2 dx\le Ed^n,
\end{equation}
and
\begin{equation}\label{E:u-normalization}
\|u\|^2_{Q_d}=\int_{Q_d}|u|^2 dx =d^n.
\end{equation}
Note that due to the diamagnetic inequality we have
\begin{equation}\label{E:mod-u-ineq}
\int_{Q_d}|\nabla |u||^2 dx\le Ed^n.
\end{equation}

For every $k\ge 0$ define a set $E_k\subset Q_d$ by
\begin{equation*}\label{E:Ek}
E_k=\{x|\,|u(x)|\ge k\},
\end{equation*}
and estimate the capacity of $E_k$.
This estimate is given in the following Lemma, 
and it can be also obtained 
from Theorem 10.1.3 in \cite{Mazya8}.

\begin{lemma}\label{L:cap-Ek} For every $k>0$
\begin{equation}\label{E:cap-Ek}
\capa(E_k)\le C_n E k^{-2} d^n. 
\end{equation}
\end{lemma}

\medskip
{\bf Proof.} Let us take $v(x)=\max(k-|u(x)|,0)$. Then $v\in\Lip(Q_d)$, 
$0\le v\le k$, and $v|_{E_k}=0$. Using Lemma \ref{L:Mazya1}, we get
\begin{equation}\label{E:cap-Ek-est}
\capa(E_k)\le \frac{C_n\int_{Q_d}|\nabla v|^2 dx}{d^{-n}\int_{Q_d}v^2 dx}.
\end{equation}
Note that $|\nabla v|\le|\nabla |u||$ almost everywhere, so \eqref{E:mod-u-ineq}
implies that
\begin{equation}\label{E:nabla-v-est}
\int_{Q_d}|\nabla v|^2 dx\le Ed^n.
\end{equation}
Let us estimate the denominator in  \eqref{E:cap-Ek-est}  from below. 
We have
\begin{equation*}
\|k\|\le \|k-|u|\|+\|u\|\le \|(k-|u|)_+\|+2\|u\|=\|v\|+2\|u\|,
\end{equation*}
where $\|\cdot\|$ is the norm in $L^2(Q_d)$.
Therefore
\begin{equation*}
\|v\|\ge \|k\|-2\|u\|=(k-2)d^{n/2}.
\end{equation*}
and the desired inequality \eqref{E:cap-Ek} follows from
\eqref{E:cap-Ek-est} and \eqref{E:nabla-v-est} provided $k\ge 3$.
It also obviously holds for $k<3$ because $E\ge d^{-2}$. $\square$ 

\bigskip
To continue the proof of Proposition \ref{P:nec}
note that the desired inequality \eqref{E:nec-ineq} holds  if 
and only if the estimate
\begin{equation}\label{E:nec-ineq-F}
\mu(Q_d;H_{a,V})\le C_nE\left(1+\frac{1}{f_n(\tmu_0)d^{n-2}}\int_{Q_d\setminus F}Vdx\right)
\end{equation}
holds for every compact $F\subset Q_d$ such that
\begin{equation}\label{E:cap-F-est}
\capa(F)\le\beta\capa(Q_d),
\end{equation}
where $\beta=cf_n(\tmu_0)$.
Let us choose  such a compact set $F$ and denote $F'=E_k\cup F$.
Then
\begin{equation}\label{E:cap-F'-beta}
\capa(F')\le \beta\capa(Q_d)+C_nEk^{-2}d^n
\end{equation}
due to the subadditivity of capacity 
and Lemma \ref{L:cap-Ek}.

We would like to apply Lemma \ref{L:cap-F-below} to the set 
$F'$.  Using \eqref{E:cap-F'-beta},
we see that it is sufficient to assume that
\begin{equation}\label{E:beta-and-k}
\beta\le c_n/2 \quad\text{and}\quad k^2\ge \frac{C_n Ed^n} {\beta\capa(Q_d)}=
\frac{\tilde C_n E d^2}{\beta},
\end{equation}
where $C_n,c_n$ are the constants from  \eqref{E:cap-F'-beta} and \eqref{E:cap-F'-cn}.
We will assume in the future that the relations \eqref{E:beta-and-k} are satisfied.
Then
\begin{equation}\label{E:cap-beta}
\capa(F')\le 2\beta\capa(Q_d).
\end{equation} 
Now we can choose a function $\psi$ as in Lemma \ref{L:cap-F-below} and define
\begin{equation}\label{E:u'}
u'=\psi u,
\end{equation} 
where $u\in \Lip(Q_d)$ satisfies \eqref{E:ha0-ineq} and
\eqref{E:u-normalization}. Clearly, $u'|_{F'}=0$
by the definition of $\psi$.

To see that we do not cut off
too much, we need to estimate the capacity of the set
\begin{equation}\label{E:R}
R=\left\{x:\,x\in Q_d, |\psi(x)|\le\frac{1}{4}\right\}.
\end{equation}
Clearly $R\supset F'$, so $\capa(R)\ge\capa(F')$. The following  Lemma
establishes an opposite estimate. 

\begin{lemma}\label{L:cap-Reps} There exists $C_n>0$ such that
\begin{equation}\label{E:cap-Reps}
\capa(R)\le C_n\capa(F').
\end{equation}
\end{lemma}

{\bf Proof.} Take $\tilde \psi=\max\{|\psi|-\frac{1}{4},0\}$,
where $\psi$ is constructed by Lemma \ref{L:cap-F-below}. Then  
$\tilde\psi|_{R}=0$, $\tilde\psi\ge 0$ and
\begin{equation}\label{E:nabla-tilde-psi-below}
\int_{Q_d}|\nabla\tilde\psi|^2dx\le \int_{Q_d}|\nabla\psi|^2dx\le 
C_n \capa(F'),
\end{equation}
where we used \eqref{E:cap-F'-psi-below}. 
On the other hand, using \eqref{E:psi2-below} we obtain
\begin{equation*}
\frac{d^n}{4}\le \int_{Q_d}|\psi|^2 dx\le \int_{Q_d}(|\tilde\psi|+\frac{1}{4})^2 dx
\le 2\int_{Q_d}|\tilde\psi|^2 dx + \frac{d^n}{8},
\end{equation*}
hence
\begin{equation*}
d^{-n}\int_{Q_d}|\tilde\psi|^2 dx\ge \frac{1}{16}
\end{equation*}
Together with \eqref{E:nabla-tilde-psi-below} and
Lemma \ref{L:Mazya1} this implies the desired inequality
\eqref{E:cap-Reps}.  $\square$

Now let us recall the following inequalities which relate the capacity 
of a compact set $F\subset Q_d$ with its Lebesgue measure $\mes F$:
\begin{equation}\label{E:cap-mes-ge-3}
\capa(F)\ge c_n[\mes F ]^{(n-2)/n}, \qquad n\ge 3,
\end{equation} 
with $c_n=\om_n^{-2/n}n^{(2-n)/n}(n-2)^{-1}$, $\om_n$ is the $(n-1)$-volume of the unit sphere
in $\R^n$;
\begin{equation}\label{E:cap-mes-2}
\capa_{\overset{\circ}Q_{d_0}}(F)\ge c_2\left[\log\frac{d_0^2}{\mes  F }\right]^{-1}, \qquad n=2,\ d_0\ge 2d,
\end{equation}
with $c_2=(4\pi)^{-1}$  (see e.g. \cite{Mazya8}, Sect. 2.2.3).
They can be rewritten as follows:
\begin{equation}\label{E:cap-mes-ge-3-alt}
\mes F \le C_n[\capa(F)]^{n/(n-2)}, \qquad n\ge 3;
\end{equation}
\begin{equation}\label{E:cap-mes-2-alt}
\mes F \le d_0^2\exp\left(-\frac{1}{C_2\capa_{\overset{\circ}Q_{d_0}}(F)}\right), \qquad n=2,\ d_0\ge 2d.
\end{equation}  

If $n=2$, then we  only need  $d_0=2d$, which will
be assumed below. Then $\capa_{\overset{\circ}Q_{d_0}}(F)=\capa_{\overset{\circ}Q_{2d}}(F)=\capa(F)$
according to our conventions. 

\begin{lemma}\label{L:int-R-estimate} Let $R$ be a compact subset in $Q_d$. If $n\ge 3$, then
\begin{equation}\label{E:int-R-estimate}
\int_R|u|^2 dx \le C_n (\mes R)^{2/n}\int_{\R^n}|\nabla u|^2 dx, \qquad u\in \Lip_c(\R^n).
\end{equation}
If $n=2$, then
\begin{equation}\label{E:int-R-estimate-2}
\int_R|u|^2 dx \le C_2 \mes R\ \log\left(\frac{4d^2}{\mes R}\right)\int_{Q_{2d}}|\nabla u|^2 dx 
\end{equation}
for any $u\in \Lip(Q_{2d})$ with $u|_{\pa Q_{2d}}=0$. (Here $Q_d$ and $Q_{2d}$ are assumed to have the same center.)
\end{lemma}

\textbf{Proof.} It is clear from the inequality $|\nabla|u||\le |\nabla u|$ 
that without loss of generality we can assume that $u\ge 0$. 
Denote for any $t\ge 0$
\begin{equation*}
N_t=\{x|\,u(x)\ge t\}\cap R.
\end{equation*}
According to Theorem 2.3.1 from \cite{Mazya8}, for any open $\Om\supset Q_d$
\begin{equation}\label{E:cap-Dir-ineq}
\int_0^\infty\capa_\Om(N_t) d(t^2)\le 4\int_\Om |\nabla u|^2 dx, \qquad u\in C_c^\infty(\Om).
\end{equation}
Using this for $\Om=\R^n$ together with  \eqref{E:cap-mes-ge-3}, we obtain for $n\ge 3$:
\begin{align*}
&\int_R u^2 dx=\int_0^\infty\mes N_t\; d(t^2)
\le (\mes R)^{2/n}\int_0^\infty (\mes N_t)^{(n-2)/n}d(t^2)\\
&\le c_n^{-1} (\mes R)^{2/n} \int_0^\infty \capa(N_t) d(t^2) 
\le 4 c_n^{-1} (\mes R)^{2/n}\int_{\R^n} |\nabla u|^2 dx,
\end{align*}
where $c_n$ is the constant from \eqref{E:cap-mes-ge-3}.  
So \eqref{E:int-R-estimate} follows with $C_n=4c_n^{-1}$. 

Let us consider the case $n=2$. We can assume $u=0$ on $\R^2\setminus Q_{2d}$.  
Using the inequalities 
\eqref{E:cap-mes-2}, \eqref{E:cap-Dir-ineq} and the fact that 
the function $\tau\mapsto \tau\log(b/\tau)$
is increasing on $(0,b/e)$, $b>0$, we obtain 
\begin{align*}
&\int_R u^2 dx=\int_0^\infty\mes N_t\; d(t^2)\\
&\le \mes R\; \log\left(\frac{4d^2}{\mes R}\right)
\int_0^\infty \left(\log\frac{4d^2}{\mes N_t}\right)^{-1}d(t^2)\\
&\le 4\pi \mes R\; \log\left(\frac{4d^2}{\mes R}\right) \int_0^\infty \capa(N_t) d(t^2)\\ 
&\le 16\pi \mes R\; \log\left(\frac{4d^2}{\mes R}\right) \int_{Q_{2d}} |\nabla u|^2 dx,
\end{align*}
so we get \eqref{E:int-R-estimate-2} with $C_2=16\pi$. $\square$

\begin{corollary}\label{C:int-R-estimate} There exist positive constants $C_n$, $n\ge 2$, such that 
if $R$ is a compact subset in $Q_d$ then for any $u\in\Lip(Q_d)$
\begin{equation}\label{E:int-R-estimate-cube}
\int_R|u|^2 dx \le C_n (\mes R)^{2/n}\left(\int_{Q_d}|\nabla u|^2 dx+d^{-2}\int_{Q_d}|u|^2 dx\right),
\end{equation}
if $n\ge 3$, and 
\begin{equation}\label{E:int-R-estimate-2-cube}
\int_R|u|^2 dx \le C_2 \mes R\ \log\left(\frac{4d^2}{\mes R}\right)
\left(\int_{Q_{d}}|\nabla u|^2 dx+ d^{-2}\int_{Q_d}|u|^2 dx\right), 
\end{equation}
if $n=2$.
 \end{corollary}

\textbf{Proof.} The result will follow if we apply Lemma \ref{L:int-R-estimate}
to the function $v=\chi U$, where $U\in \Lip(Q_{3d})$ is an extension of $u$
by reflections,
such that 
\begin{equation*}
\int_{Q_{3d}}|U|^2 dx\le 3^n \int_{Q_d} |u|^2 dx, 
\qquad \int_{Q_{3d}}|\nabla U|^2 dx\le 3^n \int_{Q_d} |\nabla u|^2 dx,   
\end{equation*}
and $\chi\in \Lip(Q_{3d})$, $\chi=1$ on $Q_d$, $\chi=0$ on  $Q_{3d}\setminus Q_{2d}$,
$0\le \chi\le 1$, $|\nabla\chi(x)|\le 2d^{-1}$  for all $x$. $\square$ 

\begin{remark} {\rm In case $n\ge 3$ another proof of the estimate \eqref{E:int-R-estimate}
can be obtained if we use the Sobolev inequality
\begin{equation}\label{E:Sobolev-ineq}
\left(\int_{\R^n}|u|^{2n/(n-2)} dx\right)^{(n-2)/n}\le C_n\int_{\R^n}|\nabla u|^2 dx,
\quad u\in\Lip_c(\R^n).
\end{equation}
(See e.g. \cite{Lieb-Loss}, Sect. 8.3.) 
By the H\"older inequality
\begin{equation*}
\int_R |u|^2 dx\le (\mes R)^{2/n}\left(\int_R |u|^{2n/(n-2)}dx\right)^{(n-2)/n}.
\end{equation*}
Combining this with \eqref{E:Sobolev-ineq}, we obtain \eqref{E:int-R-estimate}.
} 
\end{remark}

\bigskip
{\bf Proof of Proposition \ref{P:nec}.}
Let us return to the function $u$ satisfying \eqref{E:ha0-ineq} 
(hence \eqref{E:mod-u-ineq}) and \eqref{E:u-normalization}.
We would like to apply Corollary \ref{C:int-R-estimate} to the set $R$ 
defined by \eqref{E:R}
and to the function $|u|$ in order to 
establish that
\begin{equation}\label{E:int-R-above}
\int_{R}|u|^2 dx\le \frac{1}{4}\int_{Q_d}|u|^2 dx=\frac{1}{4}d^n.
\end{equation}
The inequalities in Corollary \ref{C:int-R-estimate} (applied to $|u|$)
and the diamagnetic inequality imply for this $u$ 
\begin{align}\label{E:int-Reps-cap-ge3}
&\int_{R}|u|^2 dx\le 
C_n (\mes R)^{2/n}\left(\int_{Q_d}|\nabla_a u|^2 dx+d^{-2}\int_{Q_d}|u|^2 dx\right)\\
&\le C_n(\mes R)^{2/n}(E+d^{-2})\int_{Q_d}|u|^2 dx\le 
2C_n E(\mes R)^{2/n} \int_{Q_d}|u|^2 dx, \notag 
\end{align}
if $n\ge 3$, and
\begin{equation}\label{E:int-Reps-cap-2}
\int_{R}|u|^2 dx\le 2C_2 E \mes R\; \log\left(\frac{4d^2}{\mes R}\right) \int_{Q_d}|u|^2 dx,
\end{equation} 
if $n=2$.

Note that Lemma \ref{L:cap-Reps} and \eqref{E:cap-beta} imply 
\begin{equation}\label{E:cap-beta2}
\capa(R)\le 2C_n\beta\capa(Q_d).
\end{equation}
Now for $n\ge 3$, using the estimate
\eqref{E:int-Reps-cap-ge3}, 
we see that \eqref{E:int-R-above} will follow if
\begin{equation*}
E(\mes R)^{2/n}\le \frac{c_n}{8}
\end{equation*}
with a sufficiently small $c_n>0$. Due to \eqref{E:cap-mes-ge-3-alt}, this will hold if
\begin{equation*}
E[\capa (R)]^{2/(n-2)}\le \frac{c_n}{8}
\end{equation*}
(possibly with a different $c_n$).  Recalling \eqref{E:cap-beta2}, we see that
it suffices to take
\begin{equation*}
\beta\le c_n(Ed^2)^{(2-n)/2}=c_n f_n(\mu_0 d^2)=c_n f_n(\tmu_0).
\end{equation*} 
with a small $c_n>0$. 

Now let us assume that $n=2$ and   use the estimates \eqref{E:int-Reps-cap-2}, 
\eqref{E:cap-beta2}. Taking into account that $\capa(Q_d)=\capa(Q_1)$ 
does not depend on $d$,
we see that it suffices to have
\begin{equation*}
\beta\le c_2 \left(1+\log(Ed^2)\right)^{-1}=c_2 f_2(\mu_0 d^2)=c_2 f_2(\tmu_0)
\end{equation*}
with a sufficiently small $c_2>0$. 

In both cases we see that the condition
\begin{equation}\label{E:beta-f-est}
\beta\le c_n f_n(\tmu_0)
\end{equation}
with $f_n$ as in Definition \ref{D:F},  is sufficient for the estimate 
\eqref{E:int-R-above} to hold.
Then we conclude that 
\begin{equation*}\label{E:int-Qd-Reps}
\int_{Q_d\setminus R}|u|^2 dx\ge \frac{1}{4}d^n.
\end{equation*} 
It follows that for $u'=\psi u$, as in \eqref{E:u'},
\begin{equation*}\label{E:int-u'-below}
\int_{Q_d}|u'|^2dx\ge \frac{1}{16}\int_{Q_d\setminus R_\eps}|u|^2 dx
\ge \frac{1}{64} d^n, 
\end{equation*}
whenever $\eps\in (0,1/4]$. Let us take $\eps=1/4$. Then we get
\begin{equation}\label{E:int-u'-below}
\int_{Q_d}|u'|^2dx \ge \frac{1}{64}d^n. 
\end{equation}

\medskip
Now we can use $u'$ as a test function to estimate $\mu(Q_d;H_{a,V})$. 
We obviously have
\begin{eqnarray}\label{E:mu-u'-est}
\mu(Q_d;H_{a,V})
\le \frac{h_{a,0}(u',u')_{Q_d}+(Vu',u')_{Q_d}}{\|u'\|^2_{Q_d}}\\
= \frac{\int_{Q_d}|\nabla_a u'|^2 dx+\int_{Q_d} V|u'|^2 dx}{\int_{Q_d}|u'|^2 dx}
\notag
\end{eqnarray}
Let us estimate the terms in the right hand side turn by turn.
Since $0\le\psi\le 1$, we obtain  
\begin{eqnarray*}
h_{a,0}(u',u')_{Q_d}=\int_{Q_d}|\nabla_a u'|^2 dx=
\int_{Q_d}|\psi\nabla_a u+u\nabla\psi|^2dx\\
\le  2\int_{Q_d}|\nabla_a u|^2dx+ 2\int_{Q_d}|u\nabla\psi|^2 dx.
\end{eqnarray*} 
The first term in the right hand side is estimated by $2Ed^n$ by the choice of $u$
(see \eqref{E:ha0-ineq} and \eqref{E:u-normalization}), whereas the second one
is estimated, with the use of \eqref{E:cap-F'-below}, by 
\begin{equation*}
2k^2\int_{Q_d}|\nabla\psi|^2 dx\le C_nk^2 \capa(F') d^{-n}\int_{Q_d}|\psi|^2 dx\le 
C_nk^2\capa(F').
\end{equation*}
Taking into account 
\eqref{E:cap-beta}, we see that the right hand side here is
estimated by $C_nk^2\beta\capa(Q_d)$. Now we can choose $k$ so that
the inequality \eqref{E:beta-and-k} becomes equality, i.e.
\begin{equation*}\label{E:k-choice}
k^2=\frac{\tilde C_n Ed^n}{\beta\capa(Q_d)}.
\end{equation*}
With this choice we get $k^2\capa(F')\le \tilde C_n Ed^n$, so we finally get
\begin{equation}\label{E:ha0-est-final}
h_{a,0}(u',u')_{Q_d}\le C_nEd^n.
\end{equation}

We also obviously have
\begin{eqnarray}\label{E:V-estimate}
(Vu',u')_{Q_d}=\int_{Q_d}V|u'|^2 dx\le k^2\int_{Q_d\setminus F'} Vdx\\
\le k^2\int_{Q_d\setminus F} Vdx
=\frac{\tilde C_n Ed^n}{\beta\capa(Q_d)}\int_{Q_d\setminus F}V dx,\notag
\end{eqnarray}
where we used that $V\ge 0$, $0\le \psi\le 1$ and $\psi|_{F'}=0$.

Substituting the estimates \eqref{E:ha0-est-final} and \eqref{E:V-estimate}
into \eqref{E:mu-u'-est} and taking into account \eqref{E:int-u'-below},
we obtain
\begin{equation*}\label{E:mu-HaV-est-1}
\mu(Q_d;H_{a,V})\le C_n E\left(1+\frac{1}{\beta\capa(Q_d)}\int_{Q_d\setminus F}V dx\right).
\end{equation*}
Recalling the restriction \eqref{E:beta-f-est}, we see that it is best to take
$\beta=c_n f_n(\tmu_0)$ with an appropriate (sufficiently small) constant $c_n$. 
Thus we arrive at the inequality \eqref{E:nec-ineq-F}
which proves Proposition \ref{P:nec}, hence Theorem \ref{T:discr}. $\square$

\begin{remark}\label{R:regularity}
{\rm
The condition   $ a\in L^{\infty}_{loc}(\R^n)$ 
can be substantially relaxed.
Indeed, it was only used to
guarantee that the set $\cL$ given by 
\eqref{E:cL}
(we assume that $V\ge 1$)
is precompact in $L^2(B(0,R))$ for any
$R\in(0,\infty)$.
Let us assume
that  $|a|\in M(H^1(\R^n)\to L^2_{loc}(\R^n))$,
the space  of   pointwise multipliers  mapping $H^1(\R^n)$ into
$L^2_{loc}(\R^n)$. (Here $H^1(\R^n)$ is the standard Sobolev space
of functions $u\in L^2(\R^n)$ such that $\nabla u\in L^2(\R^n)$.)
This means that for any  $R\in(0,\infty)$
$$
\int_{B(0,R)}|a|^2|v|^2dx\leq c(R)(\|\nabla v\|^2+\|v\|^2), \quad 
v\in C_c^\infty(\R^n),
$$
where $\|\cdot\|$ is the norm in $L^2(\R^n)$. 
Applying this to $v=|u|$, we obtain by the diamagnetic inequality  
$$
\int_{B(0,R)}|a|^2|u|^2dx\leq c(R) \quad {\rm for\ all} \; u\in \cL.
$$
Therefore,
$$
\|\nabla u\|^2_{L^2(B(0,R))}\leq 2\|\nabla_a u\|^2_{L^2(B(0,R))}+
2\| |a| u\|^2_{L^2(B(0,R))}\leq 2(1+ c(R)), \ u\in\cL.
$$
It remains to note that the set
$$
\{u\in C^{\infty}_c(\R^n)|\; \|\nabla u\|^2_{L^2(B(0,R))}+
\|u\|^2_{L^2(B(0,R))}\leq 3+ 2c(R)\}
$$
is precompact in   $L^2(B(0,R))$ due to the Rellich
Lemma.               
     
The  space    $M(H^1(\R^n)\to L^2_{loc}(\R^n))$ can be described
analytically in various ways 
(see   \cite{Mazya64}, Corollary 2.3.3 in \cite{Mazya8}, \cite{Kerman-Sawyer},
\cite{Mazya-Verbitsky}). 
For example, $|a|\in M(H^1(\R^n)\to L^2_{loc}(\R^n))$
if and only if for any unit ball $B(x,1)$ 
$$
\sup_F \frac{\int_F|a|^2dx}{\capa(F)}\leq c(x),
$$
where the supremum is taken over all compact  subsets $F\subset\bar B(x,1)$,
and $c=c(x)$ is continuous on $\R^n$.
 
Using the inequalities \eqref{E:cap-mes-ge-3} and \eqref{E:cap-mes-2},
we see that
it is sufficient to require that $a$ satisfies the condition
$$
\int_F|a|^2dx\leq c(x)(\mes(F))^{(n-2)/n}, \quad n>2,
$$
and
$$
\int_F|a|^2dx\leq c(x)\left( \log \frac{4}{\mes(F)}\right)^{-1},\quad n=2.
$$
It is easy to see that that the following condition on $a$
is stronger, hence also sufficient:
$a\in L^n_{loc}(\R^n)$ if $n>2$ and $|a|^2\log_+|a|\in L^1_{loc}(\R^2)$ if $n=2$.

Due to the gauge invariance it suffices that one of the conditions above
is satisfied for some $a'=a+d\phi$ with a scalar function (or a distribution)
$\phi$. 
}
\end{remark}

\section{Necessity: precision}

In this section we will construct an operator $H_{a,V}$ which will provide 
a proof of Theorem \ref{T:precision}, in particular, the precision of the exponents 
in \eqref{E:f_n}.

Let us consider a hyperplane 
\begin{equation}\label{E:hyperplane}
L=\{x|\; x^1+x^2+\dots+ x^n=0\}\subset \R^n.
\end{equation}
It divides its complement in $\R^n$ into two parts
\begin{equation}\label{E:parts}
L_\pm=\{x|\; \pm(x^1+x^2+\dots+ x^n)> 0\}.
\end{equation}
Let us take two operators $H_{\tilde a,0}$ and $H_{0,\widetilde V}$ in $\R^n$, so that each of them 
has discrete spectrum in $\R^n$, and then define $H_{a,V}$ as follows: 
\begin{equation}\label{E:HaV-example}
H_{a,V}=H_{\tilde a,0} \ {\rm in} \ L_-, \qquad H_{a,V}=H_{0, \widetilde V} \ {\rm in} \ L_+.
\end{equation} 
So $a$ and $V$ are obtained by restriction of $\tilde a$ and $\widetilde V$ to $L_-$ and $L_+$
respectively, with subsequent extensions by $0$ to the complementary half-spaces $L_+$ and $L_-$.

\ms
Theorem \ref{T:precision} will immediately follow from

\begin{proposition}\label{P:precision} 
The operator $H_{a,V}$, defined by \eqref{E:HaV-example}, has a discrete spectrum,
and satisfies the condition  formulated in Theorem \ref{T:precision}.
\end{proposition}

\textbf{Proof.} We will establish the discreteness of spectrum of $H_{a,V}$ by the 
necessary and sufficient conditions from  Theorem \ref{T:discr}.  
To this end we can use tiling cubes
with one of the faces parallel to $L$, and with interiors
in one of the half-spaces $L_{\pm}$ (see Remark \ref{R:tiling}). Then the discreteness of 
the spectrum of $H_{a,V}$ immediately follows from the corresponding properties of the operators 
$H_{\tilde a, 0}$ and $H_{0,\widetilde V}$.

Now let us choose arbitrary $d>0$, and a decreasing function $f:[0,+\infty)\to (0,1)$ satisfying
\eqref{E:not-adm-ge-3} in case $n\ge 3$ and  \eqref{E:not-adm-2} in case $n=2$. We claim then 
that the condition\eqref{E:b_f} (with $c_n=1$) is not satisfied for the cubes $Q_d$ with 
the edges parallel
to the coordinate axes (where the hyperplane $L$ has the form \eqref{E:hyperplane}).

We will consider only the cubes $Q_d$ which have ``small" intersection with 
$L_+$, with $x^1+x^2+\dots + x^n=\de >0$ at the corner of the cube where the sum
$x^1+x^2+\dots + x^n$ is maximal. We will assume that $\de\le d$. 
Then the intersection of $Q_d$ with  $\bar L_+$
(the closure of $L_+$) will be
a tetrahedron which is isometric to the tetrahedron 
\begin{equation*}
\left\{x=(x^1,\dots x^n)|\;x^j\ge 0, \sum_{j=1}^n x^j\le \de\right\}.
\end{equation*}
Clearly
\begin{equation}\label{E:cap-tetr-ge-3}
\capa(Q_d\cap \bar L_+)=c_{n}^{(1)} \de^{n-2}, \quad n\ge 3,
\end{equation} 
\begin{equation}\label{E:cap-tetr-2}
C_2^{-1} \left[\log\left(\frac{2d}{\de}\right)\right]^{-1} \le \capa(Q_d\cap \bar L_+)\le 
C_2 \left[\log\left(\frac{2d}{\de}\right)\right]^{-1},
\quad n=2.
\end{equation} 
Since $Q_d\cap \bar L_+$ is free of magnetic field ($a=0$ there) and contains a ball of diameter 
$c_n^{(2)}\de$, then, taking only test functions from $C_c^\infty(\overset{\circ}Q_d\cap L_+)$,
we obtain
\begin{equation}\label{E:mu0-delta-estimate}
\mu_0(Q_d)\le C_n\de^{-2},
\end{equation}  
if $C_n>0$ is sufficiently large.

Now we would like the sets $Q_d\cap \bar L_+$ to be negligible in the sense of 
Theorem \ref{T:discr} with the use of the function $f$, i.e.
\begin{equation}\label{E:negligible-example}
\capa(Q_d\cap \bar L_+)\le f(\mu_0 d^2) \capa(Q_d).
\end{equation}
If this is the case, then we will have $M_{\ga}(Q_d;V)=0$ and
\begin{equation}\label{E:mu0-M-estimate}
\mu_0(Q_d)+d^{-n}M_{\ga}(Q_d;V)\le C_n\de^{-2}. 
\end{equation}
with $\ga=f(\mu_0d^2)$. The condition \eqref{E:b_f} means that the left hand side of
\eqref{E:mu0-M-estimate} tends to $+\infty$ as $Q_d\to\infty$. This will not hold if we
are able to provide a sequence of cubes $Q_d\to\infty$ satisfying 
\eqref{E:mu0-M-estimate} with a fixed $\de>0$. This, in turn,  will follow if we find $\de>0$
(sufficiently small) and a sequence of cubes, constructed by the procedure above, such that 
the negligibility condition
\eqref{E:negligible-example} holds for these cubes.
 
Due to the monotonicity of $f$ and the estimate \eqref{E:mu0-delta-estimate} , 
the condition \eqref{E:negligible-example}
will follow if  we have
\begin{equation}\label{E:negligible-example2}
\capa(Q_d\cap \bar L_+)\le c_n f\left(C_n\left(\frac{\de}{d}\right)^{-2}\right) d^{n-2},
\end{equation}
where $c_n=\capa(Q_1)$. Now using \eqref{E:cap-tetr-ge-3} and 
\eqref{E:not-adm-ge-3} in case $n\ge 3$ we can rewrite this
condition in the form
\begin{equation}\label{E:E:negligible-example3}
c_n^{(1)}\left(\frac{\de}{d}\right)^{n-2}\le 
c_n\left(1+C_n\left(\frac{\de}{d}\right)^{-2}\right)^{(2-n)/2} 
h\left(C_n\left(\frac{\de}{d}\right)^{-2}\right),
\end{equation}
so it obviously holds if $\de/d$ is sufficiently small, because $h(t)\to +\infty$ as $t\to +\infty$.

In case $n=2$, due to \eqref{E:cap-tetr-2} and \eqref{E:not-adm-2}, 
the inequality \eqref{E:negligible-example2} will be
fulfilled if we require that
\begin{equation}\label{E:negligible-example4}
C_2\left[\log\left(\frac{2d}{\de}\right)\right]^{-1}\le 
\left[1+\log\left(C_2\left(\frac{\de}{d}\right)^{-2}\right)\right]^{-1} 
h\left(C_2^{-1}\left(\frac{\de}{d}\right)^{-2}\right)
\end{equation} 
for a sufficiently large $C_2>0$. This again holds if $\de/d$ is sufficiently small. $\square$

\section{Positivity}

In this section we will prove Theorem \ref{T:pos}. We will consider operators $H_{a,V}$ with
$V\in L^1_{loc}(\R^n)$, $V\ge 0$, $a\in L^2_{loc}(\R^n)$.

The proof will be essentially based
on the same arguments as the proof of Theorem \ref{T:discr}, except that the large cubes are essential
here (instead of small cubes). 

We will use the notations 
from Section \ref{S:Intro} and start with the following localization result:

\begin{proposition}\label{P:pos-loc}
For an operator $H_{a,V}$ the following conditions are equivalent:

\medskip
$(a)$ There exists $d_1>0$ such that $H_{a,V}\ge d_1^{-2}I$, or, equivalently, 
$0$ is not in the  spectrum of $H_{a,V}$ in $L^2(\R^n)$ 
(i.e. the spectrum is in $[\eps_0,+\infty)$ for some $\eps_0>0$).

\medskip
$(b)$ There exist $d>0$ and $d_1>0$ such that $\mu(Q_d;H_{a,V})\ge d_1^{-2}$ 
for every cube $Q_d\subset\R^n$.

\medskip
$(c)$ There exist $d_1>0$ and $d_2>0$ such that 
for every $d>d_2$ we have $\mu(Q_d;H_{a,V})\ge d_1^{-2}$ 
for every cube $Q_d\subset\R^n$.

\medskip
$(d)$ There exists $d_1>0$ such that for every $d>0$ we have $\la(Q_d;H_{a,V})\ge d_1^{-2}$
for every cube $Q_d\subset\R^n$.

\medskip
$(e)$ There exists $d_1>0, d_2>0$ such that for every $d>d_2$ we have 

\noindent
$\la(Q_d;H_{a,V})\ge d_1^{-2}$
for every cube $Q_d\subset\R^n$.
\end{proposition}

{\bf Proof.} The equivalence of $(a)$, $(d)$ and $(e)$ follows from the fact that the
quadratic form $h_{a,V}$ of $H_{a,V}$ is obtained as the closure 
from the original domain $C_c^\infty(\R^n)$. 

Using the inequality (\cite{Molchanov, Kondrat'ev-Shubin, Kondratiev-Shubin})
\begin{equation}\label{E:la-mu}
\mu(Q_d;H_{a,V})\le \la(Q_d;H_{a,V})\le A_n \mu(Q_d;H_{a,V})+\frac{B_n}{d^2},
\end{equation}
where $A_n>0$, $B_n>0$, we immediately see that $(d)$ implies that
\begin{equation*}
\mu(Q_d;H_{a,V})\ge A_n^{-1}\left[\la(Q_d;H_{a,V})-\frac{B_n}{d^{2}}\right]\ge 
A_n^{-1}\left[\frac{1}{d_1^2}-\frac{B_n}{d^{2}}\right]
\ge \frac{1}{2A_n d_1^2}, 
\end{equation*}
if $d>d_2>0$ with $d_2^2\ge 2B_n d_1^2$. So $(d)$ implies  $(c)$. 
Obviously $(c)$ implies $(b)$.

Now we see that the Proposition will be proved if we establish
that $(b)$ implies $(a)$.
So let us assume that $(b)$ holds. Then we have
\begin{equation}\label{E:pos-Qd}
\|u\|^2_{Q_d}\le d_1^2 h_{a,V}(u,u)_{Q_d}, \quad u\in\Lip(Q_d),
\end{equation}
for every cube $Q_d$ with $d>0$ taken from the condition $(b)$.
If we take an arbitrary $u\in \Lip_c(\R^n)$ and sum up the inequalities \eqref{E:pos-Qd}  over a tiling of
$\R^n$ by cubes $Q_d$, we will get the inequality
$\|u\|^2\le d_1^2 h_{a,V}(u,u)$ which proves $(a)$. $\square$

\bigskip
{\bf Proof of Theorem \ref{T:pos}.} Clearly the following implications hold:

\noindent
$$
({\rm e})\Longrightarrow({\rm c})\Longrightarrow ({\rm b})\qquad
{\rm and}
\qquad
({\rm e})\Longrightarrow({\rm d})\Longrightarrow ({\rm b}).
$$

So it suffices to prove the following two implications:

$({\rm b})\Longrightarrow({\rm a})$ (sufficiency of $({\rm b})$) and
$({\rm a})\Longrightarrow({\rm e})$ (necessity of $({\rm e})$).

\bigskip
{\bf Proof of the implication $({\rm b})\Longrightarrow({\rm a})$.} Let us assume that there exist
$c>0$, $d_1>0$ and $d>0$ 
such that the inequality
\eqref{E:M-pos} holds for all cubes $Q_d$.

The desired strict positivity will follow if we prove the inequality
\begin{equation}\label{E:pos-ineq}
\int_{\R^n}|u|^2dx\le d_2^2 \int_{\R^n}\left(|\nabla_au|^2+V|u|^2\right)dx, \quad
u\in C_c^\infty(\R^n).
\end{equation}
Note first that  for every $u\in\Lip(Q_d)$
\begin{equation}\label{E:pos-mag}
\mu_0(Q_d)\int_{Q_d}|u|^2 dx\le \int_{Q_d}|\nabla_a u|^2 dx\le 
\int_{Q_d}\left(|\nabla_au|^2+V|u|^2\right)dx.
\end{equation}

As we did in the proof of Proposition \ref{P:general-suff}, let us split the cubes $Q_d$
from a tiling of $\R^n$ into two types:

Type I ({\it large energy of the magnetic field in $Q_d$}):
\begin{equation*}
\mu_0(Q_d)> \frac{1}{2d_1^2};
\end{equation*}

Type II ({\it small energy of the magnetic field in $Q_d$}):

\begin{equation*}
\mu_0(Q_d)\le \frac{1}{2d_1^2}.
\end{equation*}

For a type I cube $Q_d$ we obtain from \eqref{E:pos-mag} that for every $u\in\Lip(Q_d)$
the inequality \eqref{E:pos-Qd} holds with $2d_1^2$ instead of $d_1^2$.

Now let $Q_d$ be a type II cube. Then we have 
\begin{equation*}
d^{-n}M_c(Q_d;V)\ge \frac{1}{2d_1^2}. 
\end{equation*} 
Due to Lemma \ref{L:Mazya2} and Remark \ref{R:main-ineq} we obtain
for every $u\in \Lip(Q_d)$ and  $c>0$
\begin{equation*}\label{E:Mazya2-mag}
\int_{Q_d}|u|^2 dx\le {C_nd^2\over c}\  \int_{Q_d}|\nabla_a u|^2 dx
+{4d^n\over M_c(Q_d;V)}\
\int_{Q_d}|u|^2 Vdx,
\end{equation*}
and  we get
\begin{equation*}\label{E:Mazya2-mag}
\int_{Q_d}|u|^2 dx\le Cd^2 \int_{Q_d}|\nabla_a u|^2 dx
+8d_1^2 \int_{Q_d}|u|^2 Vdx,
\end{equation*}
where $C=C_n/c$. Taking $d_2>0$ such that
\begin{equation*}
d_2^2=\max\left(Cd^2, 8d_1^2\right),
\end{equation*}
we obtain \eqref{E:pos-Qd} with $d_2^2$ instead of $d_1^2$.

So we obtained the inequalities \eqref{E:pos-Qd} (with $d_2^2$ instead of $d_1^2$)
for both types of cubes. This means that the condition $(b)$ in Proposition \ref{P:pos-loc}
is satisfied, hence the spectrum of $H_{a,V}$ is discrete.  $\square$

\bigskip
{\bf Proof of the implication $({\rm a})\Longrightarrow({\rm e})$ .}  We will use Proposition \ref{P:nec}
in the same way as in the proof of Theorem \ref{T:discr}. Recall the  notation
$E=\mu_0(Q_d)+d^{-2}$ which was introduced in the formulation of Proposition \ref{P:nec},
and will be used here too, though for large $d$ when the difference between 
$E$ and $\mu_0(Q_d)$ becomes small.

According to Proposition \ref{P:pos-loc} we can assume that its condition $(c)$ is satisfied,
i.e. $\mu(Q_d;H_{a,V})\ge d_3^{-2}$ for every cube $Q_d$ with $d>d_4$, where
$d_3,d_4>0$ are sufficiently large. Then 
due to \eqref{E:nec-ineq} we have for such $d$
\begin{equation}\label{E:mu-plus-M}
\mu_0(Q_d)+\frac{E}{f(\tmu_0)d^{n-2}}M_{c_nf_n(\tmu_0)}(Q_d;V)\ge 
\frac{1}{C_nd_3^{2}}-\frac{1}{d^2}\ge \frac{1}{d^2},
\end{equation}
provided $d^2\ge 2C_n d_3^2$.

Now note that in the case when  
\begin{equation*}
\mu_0(Q_d)\ge \frac{1}{d^2},
\end{equation*}
the desired inequality \eqref{E:M-pos-dn} becomes obvious
(with $\tc_n=1$). 
So from now on we can assume that
\begin{equation*}
\mu_0(Q_d)\le \frac{1}{d^2}.
\end{equation*}
This implies that 
\begin{equation*}
f_n(\tmu_0)=f_n(\mu_0 d^2)\ge f_n(1)>0, \qquad n\ge 2.
\end{equation*}
We also have in this case $E\le 2d^{-2}$.
It follows that the coefficient in front of $M_{c_nf_n(\tmu_0)}(Q_d;V)$ in 
\eqref{E:mu-plus-M} is bounded from above by $C_nd^{-n}$.
Hence the left hand side in \eqref{E:mu-plus-M} is bounded from above by
$\tilde C_n[\mu_0(Q_d)+d^{-n}M_{c_n}(Q_d;V)]$ and the desired inequality
\eqref{E:M-pos-dn} follows with $\tc_n=\min(\tilde C_n^{-1},1)$.
This ends the proof of Theorem \ref{T:pos}. $\square$

\section{Operators in domains}\label{S:domains}

In this section we will discuss the discreteness of spectrum and 
strict positivity for the magnetic Schr\"odinger operators in arbitrary open 
subsets $\Omega\subset\R^n$ with the Dirichlet boundary conditions on $\pa\Omega$.
It occurs that the methods developed above can be extended 
to this case and provide necessary and sufficient conditions 
so that the results of the previous sections appear as a particular case
when $\Omega=\R^n$. 
Note that the geometry of the domain may
contribute to the discreteness of spectrum or strict positivity
and even be the only cause of these properties.

Let $H_{a,V}$ be the magnetic Schr\"odinger operator defined as in Section \ref{S:Intro}
but in $L^2(\Omega)$. We will assume that $V\in L^1_{loc}(\Omega)$, $V\ge 0$, $a\in L^2_{loc}(\Om)$.
For the discreteness of spectrum results we will assume that $a$ is bounded in $\Omega\cap B(0,R)$ for every
$R>0$, though this condition may be substantially weakened as explained in Remark \ref{R:regularity}. The operator
$H_{a,V}$ is defined by the quadratic form
\eqref{E:haV} on functions
$u\in C_c^\infty(\Om)$.

We will define the Molchanov functional in $\Om$ as follows
\begin{equation*}
M_{\ga,\Om}(Q_d;V)=
\inf_F\left\{\int_{Q_d\setminus F}Vdx| \capa(F)\le\ga\capa(Q_d), 
F\supset Q_d\cap(\R^n\setminus\Om)\right\},\notag
\end{equation*}
where $0<\ga<1$, $F$ is a closed subset in $Q_d$. By definition it is $+\infty$ 
if there is no
sets $F$ satisfying the condition in the braces, i.e. if 
\begin{equation}\label{E:noF}
\capa(Q_d\cap(\R^n\setminus\Om)) > \ga\capa(Q_d).
\end{equation}
The numbers $\la(Q_d;H_{a,V})$ and $\mu(Q_d;H_{a,V})$ should be replaced by 
the numbers $\la_\Om(Q_d;H_{a,V})$ and $\mu_\Om(Q_d;H_{a,V})$ which are defined by
the same formulas \eqref{E:lambda}, \eqref{E:mu} (with $G=Q_d$) but with an additional requirement
on $u$ to vanish in a neighborhood of $Q_d\cap(\R^n\setminus \Om)$. 
Then the same localization results (see e.g. Theorems 1.1--1.3 in \cite{Kondratiev-Shubin}) hold.
For example,  $H_{a,V}$ has a discrete spectrum in $L^2(\Om)$ if and only if for any fixed $d>0$
\begin{equation*}
\mu_\Om(Q_d;H_{a,V})\to +\infty \quad {\rm as}\quad Q_d\to\infty.
\end{equation*} 

The appropriate modification of $\mu_0$ (the local energy of the magnetic field) is
\begin{equation*}
\mu_{0,\Om}=\mu_{0,\Om}(Q_d)=\mu_{0,\Om}(Q_d;a)=\mu_\Om(Q_d;H_{a,0}).
\end{equation*} 

With these notations the following theorems are obtained by simple repetition
of arguments given in the previous sections.

\begin{theorem}\label{T:discr-dom} 
Theorem \ref{T:discr} holds for $H_{a,V}$ in $L^2(\Om)$ if we replace 
$\mu_0$ by $\mu_{0,\Om}$ and $M_\gamma(Q_d;V)$ by $M_{\ga,\Om}(Q_d;V)$.
\end{theorem}  

\begin{theorem}\label{T:pos-dom}
Theorem \ref{T:pos} holds for $H_{a,V}$ in $L^2(\Om)$ if we replace 
$\mu_0$ by $\mu_{0,\Om}$ and $M_\gamma(Q_d;V)$ by $M_{\ga,\Om}(Q_d;V)$.
\end{theorem}

The appropriate modifications of Corollaries \ref{C:necessary} and \ref{C:sufficient}
hold as well. The same replacements of $\mu_0$ by $\mu_{0,\Om}$ and $M_c$ by $M_{c,\Om}$ 
should be made in the formulations, and the integral in
\eqref{E:cor-nec} should be replaced by $M_{0,\Om}(Q_d;V)$ which is equal to this integral if 
$Q_d\subset \Om$ and to $+\infty$ otherwise.

Now we will formulate some more specific corollaries of Theorem \ref{T:discr-dom}, 
which treat the cases when one or both fields vanish.
We will start with the case when $a\equiv 0, V\equiv 0$. 

\begin{corollary}\label{C:noaV}
There exists $c_n>0$ such that for every function $g\in\cG$ (see Definition \ref{D:F}) 
the following conditions are equivalent:

$({\rm a})$ The spectrum of the operator $H_{0,0}=-\De$ in $L^2(\Om)$ with the Dirichlet boundary
conditions on $\pa\Om$ is discrete.

$({\rm b}_g)$ $\exists$ $d_0>0$,  $\forall$ $d\in(0,d_0)$, $\exists$ $R=R(d)>0$,
$\forall$ $Q_d$ such that 

\noindent
$Q_d\cap(\R^n\setminus B(0,R))\ne\emptyset$, the inequality \eqref{E:noF} is satisfied
with $\ga=c_n g(d)^{-1}d^2$.

In particular, all conditions $({\rm b}_g)$ for different $g\in\cG$ are equivalent.
\end{corollary}

Instead of $({\rm b}_g)$ we can equivalently write that $\exists$ $d_0>0$,  $\forall$ $d\in(0,d_0)$
\begin{equation*}
\underset{Q_d\to\infty}{\lim\inf}\;\frac{\capa(Q_d\cap (\R^n\setminus \Om))}{\capa(Q_d)}>\ga,
\end{equation*}
with the same $\ga$ as above (we can replace $c_n$ by a smaller positive number).

Note that the condition $({\rm b}_g)$ is a purely geometric condition on the open set $\Om\subset\R^n$.
The equivalence of these conditions for different functions $g\in\cG$ is a non-trivial 
geometric property of the capacity.

The next corollary treats the case when $a\equiv 0$, i.e. there is no magnetic field.

\begin{corollary}\label{C:noa}
There exists $c_n>0$ such that for every  $g\in\cG$ 
the following conditions are equivalent:

$({\rm a})$ The spectrum of the operator $H_{0,V}=-\De+V$ in $L^2(\Om)$ with the Dirichlet boundary
conditions on $\pa\Om$ is discrete.

$({\rm b}_g)$ $\exists$ $d_0>0$,  $\forall$ $d\in(0,d_0)$
\begin{equation*}
M_{\ga,\Om}(Q_d;V)\to +\infty \quad {\rm as} \quad Q_d\to\infty,
\end{equation*}
where $\ga=c_n g(d)^{-1}d^2$.

$({\rm c}_g)$ $\exists$ $d_0>0$,  $\forall$ $d\in(0,d_0)$
\begin{equation*}
\underset{Q_d\to\infty}{\lim\inf}\; d^{-n} M_{\ga,\Om}(Q_d;V)\ge g(d)^{-1},
\end{equation*}
with the same $\ga$ as in $({\rm b}_g)$.

In particular, all conditions $({\rm b}_g), ({\rm c}_g)$ for different $g\in\cG$ 
are equivalent.
\end{corollary}

Finally, we consider the case when $V\equiv 0$.
To this end we need the quantity $\mu_{0,\Om}^{(\ga)}(Q_d)$ which is defined as
$\mu_{0,\Om}(Q_d)$ if $\capa(Q_d\cap(\R^n\setminus\Om))\le\ga\capa(Q_d)$
and $+\infty$ otherwise (i.e. if \eqref{E:noF} is satisfied). 
 
\begin{corollary}\label{C:noV}
There exists $c_n>0$ such that for every  $n$-admissible pair $(f,g)$ 
(see Definition \ref{D:F})
the following conditions are equivalent:

$({\rm \tilde a})$ The spectrum of the operator $H_{a,0}$ in $L^2(\Om)$ with the Dirichlet boundary
conditions on $\pa\Om$ is discrete.

$({\rm \tilde b}_g)$ $\exists$ $d_0>0$,  $\forall$ $d\in(0,d_0)$
\begin{equation*}
\mu_{0,\Om}^{(\ga)}(Q_d)\to +\infty \quad {\rm as} \quad Q_d\to\infty,
\end{equation*}
where $\ga=c_n f(\mu_{0,\Om}d^2)g(d)^{-1}d^2$.

$({\rm \tilde c}_g)$ $\exists$ $d_0>0$,  $\forall$ $d\in(0,d_0)$
\begin{equation*}
\underset{Q_d\to\infty}{\lim\inf} \mu_{0,\Om}^{(\ga)}(Q_d)\ge g(d)^{-1},
\end{equation*}
with the same $\ga$ as in $({\rm \tilde b}_g)$.

In particular, all conditions $({\rm \tilde b}_g), ({\rm \tilde c}_g)$ 
for different $g\in\cG$  are equivalent.
\end{corollary}

We skip formulations of similar Corollaries of Theorem \ref{T:pos-dom}.

\section*{Appendix: Proofs of Lemmas \ref{L:Mazya1}, \ref{L:Mazya2}, \ref{L:cap-F-below}.}\label{S:appendix}
\renewcommand{\thesection}{A}
\setcounter{theorem}{0}
\setcounter{equation}{0}

In this appendix, for the convenience of the readers, we will provide proofs
of Lemmas \ref{L:Mazya1}, \ref{L:Mazya2}, \ref{L:cap-F-below}. These proofs are
simpler compared with the proofs given in \cite{Mazya8} due to the fact that
the corresponding results in \cite{Mazya8} have much bigger generality.  

\ms
Let us recall the following classical {\it Poincar\'e inequality}
(see e.g. \cite{Gilbarg-Trudinger}, Sect. 7.8, or \cite{Kondrat'ev-Shubin}, Lemma 5.1):
\begin{equation}\label{E:Poincare}
||u -\bar u||^2_{Q_d}
\le {d^2\over \pi^2} \int_{Q_d} |\nabla u(x)|^2 dx,
\end{equation}
where $\|\cdot\|_{Q_d}$ is the norm in $L^2(Q_d)$, $u\in \Lip(Q_d)$, and
\begin{equation*}\label{E:mean-def}
\bar u = d^{-n} \int_{Q_d} u(x)\, dx
\end{equation*}
is the mean value of $u$ on $Q_d$.
 
\bs
{\bf Proof of Lemma \ref{L:Mazya1}.}
Let us normalize $u$ by 
\begin{equation*}\label{E:normalization}
d^{-n}\int_{Q_d} |u(x)|^2 dx=1,
\end{equation*}
i.e. $\overline {|u|^2}=1$ (we will call it {\it the standard normalization}). 
By the Cauchy-Schwarz inequality we obtain
\begin{equation}\label{E:|u|}
\overline{|u|}\le \left(\overline{|u|^2}\right)^{1/2}=1
\end{equation} 

Replacing $u$ by $|u|$ does not change the denominator
and may only decrease the numerator in \eqref{E:cap-above}. Therefore
we can restrict ourselves to Lipschitz functions $u\ge 0$.  

Let us denote $\phi=1-u$. Then $\phi=1$ on $F$, and $\bar\phi=1-\bar u\ge 0$ due to \eqref{E:|u|}.
Let us estimate $\bar\phi$ from above. 
Obviously
\begin{equation*}\label{E:1-|u|}
\bar\phi=d^{-n/2}(\|u\|-\|\bar u\|)
\le d^{-n/2}\|u-\bar u\|,
\end{equation*}
where $\|\cdot\|$ means the norm in $L^2(Q_d).$
So the Poincar\'e inequality gives
\begin{equation*}\label{E:1-bar-u}
\bar\phi
\le \pi^{-1} d^{-n/2+1}\|\nabla u\|=\pi^{-1} d^{-n/2+1}\|\nabla \phi\|, 
\end{equation*}
hence
\begin{equation*}\label{E:phi-square}
\bar{\phi}^2
\le {1\over \pi^2} d^{2-n} \int_{Q_d} |\nabla\phi|^2dx.
\end{equation*}

Using the Poincar\'e inequality again, we obtain
\begin{equation*}\label{E:phi-square-est}
\|\phi\|^2=\|(\phi-\bar\phi)+\bar\phi\|^2\le
2\|\phi-\bar\phi\|^2+2\|\bar\phi\|^2\le
{4d^2\over \pi^2} \int_{Q_d} |\nabla\phi |^2dx,
\end{equation*}
or
\begin{equation}\label{E:Mazya(1)}
\int_{Q_d} \phi^2 dx\le {4d^2\over \pi^2} \int_{Q_d}
|\nabla \phi |^2 dx. 
\end{equation}
Let us extend $\phi$ outside $Q_d$
by symmetries in the faces of $Q_d$, so that
the extension $\tilde \phi$ satisfies
\begin{equation*}\label{E:Q3d}
\int_{Q_{3d}} |\nabla \tilde\phi |^2 dx=3^n\int_{Q_d}
|\nabla \phi|^2dx,\  \int_{Q_{3d}}|\tilde\phi|^2dx=3^n \int_{Q_d}
|\phi|^2dx.
\end{equation*}
Denote by $\eta$ a continuous piecewise linear function, such that
$\eta=1$ on $Q_d$, $\eta=0$ outside $Q_{2d}$,
$0\le \eta\le 1$ and $|\nabla\eta|\le 2d^{-1}$.
Then
\begin{equation*}\label{E:capF-est}
\capa(F)\le 
\int_{Q_{2d}}|\nabla (\tilde\phi\eta)|^2dx
\le 2\cdot 3^n\left(\int_{Q_d}|\nabla \phi|^2dx+
4d^{-2}\int_{Q_d} \phi^2 dx\right).
\end{equation*}
Taking into account that $|\nabla\phi|=|\nabla u|$
and using \eqref{E:Mazya(1)}, we obtain
\begin{equation*}
\capa(F)\le C_n \int_{Q_d} |\nabla u|^2 dx,
\end{equation*}
which is equivalent to the desired estimate \eqref{E:cap-above}. $\square$

\bs
{\bf Proof of Lemma \ref{L:Mazya2}.} Let ${\cal M}_\tau=\{x\in Q_d:|u(x)|> \tau\},$
where $\tau\ge 0$.
Since
\begin{equation*}
|u|^2\le 2\tau^2+2(|u|-\tau)^2 \quad \hbox{on}\ {\cal M}_\tau,
\end{equation*}
we have for all $\tau$
\begin{equation*}
\int_{Q_d} |u|^2dx\le 2\tau^2d^n+
2\int_{\cal M_\tau} (|u|-\tau)^2dx.
\end{equation*}
Let us take 
\begin{equation*}
\tau^2={1\over 4d^n} \int_{Q_d} |u|^2 dx,
\end{equation*} 
i.e.
$\tau={1\over 2}\left(\overline{|u|^2}\right)^{1/2}$. Then for
this particular value of $\tau$ we obtain
\begin{equation}\label{E:(3)}
\int_{Q_d} |u|^2 dx\le
4\int_{\cal M_\tau} (|u|-\tau)^2dx. 
\end{equation}
Assume first that
$\capa(Q_d\setminus {\cal M}_\tau)\ge\gamma \capa(Q_d)$.
Using \eqref{E:(3)} and applying Lemma \ref{L:Mazya1} to the function $(|u|-\tau)_+$, which
equals $|u|-\tau$ on ${\cal M}_\tau$ and $0$ on
$Q_d\setminus {\cal M}_\tau$, we see that
\begin{equation*}
\capa(Q_d\setminus{\cal M}_\tau)\le 
\frac{C_n\int_{{\cal M}_\tau}|\nabla(|u|-\tau)|^2 dx}
{d^{-n}\int_{Q_d}|u|^2 dx}\le
\frac{C_n\int_{Q_d}|\nabla u|^2 dx}
{d^{-n}\int_{Q_d}|u|^2 dx}\;,
\end{equation*}
where $C_n$ is 4 times the one in \eqref{E:cap-above}. Therefore
\begin{equation*}\label{E:(4)}
\int_{Q_d}|u|^2 dx \le\frac{C_n d^n\int_{Q_d}|\nabla u|^2 dx}{\capa(Q_d\setminus{\cal M}_\tau)}
\le\frac{C_n d^n\int_{Q_d}|\nabla u|^2 dx}{\gamma\,\capa(Q_d)}
\end{equation*}
Taking into account that  $\capa(Q_d)=c_n d^{n-2}$ we see that
\begin{equation}\label{E:(4)}
\int_{Q_d}|u|^2 dx 
\le\frac{C_n d^2}{\gamma}\int_{Q_d}|\nabla u|^2 dx
\end{equation}
with yet another constant $C_n$.

Now consider the opposite case
$\capa(Q_d\setminus {\cal M}_\tau)\le \gamma\, \capa(Q_d)$.
Then we can write
\begin{align*}
&\int_{Q_d} |u|^2 V dx\ge \int_{{\cal M}_\tau} |u|^2 V dx\ge
\tau^2\int_{{\cal M}_\tau} V dx
={1\over 4d^n}\int_{Q_d}|u|^2dx\cdot \int_{{\cal M}_\tau} V dx\\
&\ge {1\over 4d^n}\int_{Q_d}|u|^2dx\cdot
\inf_F\int_{Q_d\setminus F} Vdx,
\end{align*}
where the infimum should be taken over all compact sets $F\subset Q_d$
such that $\capa(F)\le\gamma\,\capa(Q_d)$, so it becomes $M_\gamma(Q_d;V)$.
Finally we obtain in this case
\begin{equation}\label{E:(5)}
\int_{Q_d} |u|^2 dx\le 
\frac{4 d^n}{M_\ga(Q_d;V)}\int_{Q_d} V|u|^2 dx.
\end{equation}

The resulting inequality \eqref{E:Mazya(2)} follows from \eqref{E:(4)} and \eqref{E:(5)}. $\square$

\bs
{\bf Proof of Lemma \ref{L:cap-F-below}.} 
We start with a function $\phi\in \Lip_c(\R^n)$
such that $0\le \phi\le 1$, 
$\phi=1$ in a neighborhood of $F'$,
$\phi=0$ outside $Q_{d_0}$
(where for $n=2$ we take $d_0=2d$),
 and
\begin{equation}\label{E:cap-F'-phi-below}
\capa(F')\ge c_n'\int_{Q_{d_0}}|\nabla \phi|^2dx
\end{equation}
with $c_n'>0$.
It follows that
\begin{equation*}\label{E:cap-F'-phi-below-Qd}
\capa(F')\ge c_n'\int_{Q_{d}}|\nabla \phi|^2dx.
\end{equation*}
Now take $\psi=1-\phi$, so $0\le \psi\le 1$ and $\psi|_{F'}=0$. 
Then $|\nabla\psi|=|\nabla \phi|$, hence the condition
\eqref{E:cap-F'-psi-below} is obviously satisfied.
Now our goal  will be achieved if we prove that 
\eqref{E:psi2-below} holds
provided \eqref{E:cap-F'-cn} is satisfied with a sufficiently 
small  $c_n>0$. 

To prove \eqref{E:psi2-below}, note first that Lemma \ref{L:int-R-estimate} with $R=Q_d$ gives
\begin{equation}\label{E:Hardy-type}
\int_{Q_d}|\phi|^2dx\le C_n d^2\int_{Q_{d_0}}|\nabla \phi|^2dx,
\end{equation}

Hence, using  
\eqref{E:cap-F'-phi-below}, we obtain
\begin{equation*}
\overline{\phi^2}=d^{-n}\int_{Q_d}\phi^2 dx\le C_n d^{2-n}\int_{Q_{d_0}}|\nabla\phi|^2dx
\le  \frac{\tilde C_n (c_n')^{-1}\capa(F')}{\capa(Q_d)}\le \tilde C_n (c_n')^{-1}c_n,
\end{equation*}
where $c_n$ is the constant from \eqref{E:cap-F'-cn}. Now we can adjust $c_n$
so that we have $\tilde C_n (c_n')^{-1}c_n\le 1/4$. Then \eqref{E:psi2-below} follows  
from the triangle inequality. $\square$

\end{document}